\documentstyle[amsfonts,amssymb]{article} 
\newtheorem{theorem}{Theorem}[section]
\newtheorem{corollary}[theorem]{Corollary}
\newtheorem{lemma}[theorem]{Lemma}
\newtheorem{proposition}[theorem]{Proposition}

\newtheorem{definition}[theorem]{Definition}
\newtheorem{remark}[theorem]{Remark}

\newcommand{\wt}{\widetilde}

\newcommand{\mb}{\mathbb }

\newcommand{\dsp}{\displaystyle} 
\newcommand{\pa}{\partial} 
\newcommand{\be}{\begin{equation}}
\newcommand{\ee}{\end{equation}}
\newcommand{\beq}{\begin{eqnarray}}
\newcommand{\eeq}{\end{eqnarray}}
\newcommand{\beqst}{\begin{eqnarray*}}
\newcommand{\eeqst}{\end{eqnarray*}}

\date{}

\title{On the global solutions of  the Higgs boson  equation}

\author{Karen Yagdjian}

\begin{document}
\maketitle
\thispagestyle{empty}

 \centerline{Department of Mathematics, 
University of Texas-Pan American,}
   \centerline{1201 W.~University Drive, 
Edinburg, TX 78541-2999,  
USA }
 \centerline{yagdjian@utpa.edu}

\bigskip

\begin{abstract}
  In this article we study      global in time (not necessarily small)  solutions of the equation for the Higgs boson   in the Minkowski and 
  in the de~Sitter spacetimes. We reveal some qualitative behavior of the global solutions. In particular,  we formulate 
  sufficient conditions for the existence of the zeros  of  global solutions in the interior of their supports, and, 
consequently, for the creation of the so-called bubbles, which have been studied in  particle physics and inflationary  cosmology. 
We also give some  sufficient conditions for the global solution to be an oscillatory in time solution. 
\end{abstract}

\section{Introduction}

  In this article we study  the global in time,  not necessarily small, solutions of the  Higgs  boson   equation    in the Minkowski and in the 
  de~Sitter spacetime.  
The Higgs  boson  plays a fundamental role in unified theories of  weak, strong, and electromagnetic
interactions \cite{Weinberg}.

In the model  of the universe proposed by  de~Sitter,  the line
element  has the   
form
\[
ds^2=-\left( 1-\frac{2M_{bh}}{r} - \frac{\Lambda r^2}{3}\right)c^2dt^2 + \left( 1-\frac{2M_{bh}}{r} - \frac{\Lambda r^2}{3}\right)^{-1}dr^2 + r^2(d \theta ^2+ \sin^2 d\phi ^2).
\]
The constant  $M_{bh}$ may have a meaning of the ``mass of the black hole''. 
The corresponding metric with this line element 
is called the Schwarzschild - de~Sitter metric. 
In the present paper we focus on the  limit case; namely, we  set \,$M_{bh}=0$ \,to ignore completely any influence of the black hole.  
Thus,  the  line element in the de~Sitter spacetime has  the form
\[
ds^2= - \left( 1- \frac{r^2}{R^2}\right) c^2\, dt^2+ \left( 1- \frac{r^2}{R^2}\right)^{-1}dr^2 + r^2(d\theta ^2 + \sin^2 \theta \, d\phi ^2)\,.
\]
The  Lama{\^i}tre-Robertson transformation \cite{Moller} 
leads to the following form for the line element \cite[Sec.134]{Moller}: 
\[
ds^2= -   c^2\, d{t'}^2+ e^{2ct'/R}( d{x'}^2   + d{y'} ^2 +  d{z '}^2)\,.
\]
Here $R$ is the ``radius'' of the universe.   
In the Robertson-Walker spacetime \cite{Hawking}, one can choose coordinates so that the metric has the form
\[
ds^2=-dt^2+S^2(t)d \sigma ^2\,.
\]
In particular, the  metric in   
de Sitter   spacetime in the Lama{\^i}tre-Robertson coordinates \cite{Moller} has this form with the cosmic scale factor $S (t)=e^{t} $.

\smallskip

The matter waves in the de~Sitter spacetime are described by the function  $\phi $, which  satisfies equations of 
motion.  
In the  de~Sitter universe  the equation for the scalar field with  potential function \, $V$   \, 
is the covariant wave equation  
\[
\square_g \phi    = V'(\phi ) \quad \mbox{\rm or} \quad \frac{1}{\sqrt{|g|}}\frac{\partial }{\partial x^i}
\left( \sqrt{|g|} g^{ik} \frac{\partial \phi  }{\partial x^k} \right)  =V'(\phi ) \,,
\]
with the usual summation convention. Written explicitly in coordinates  in the  de Sitter spacetime it, in particular, 
for \,
\[
V'(\phi )=  -\mu^2  \phi + \lambda  |\phi |^{p-1}\phi,\qquad p>1, 
\]
 has  the form  
\begin{equation}
\label{1.1}
  \phi_{tt} +   n   \phi_t - e^{-2 t} \Delta  \phi =    \mu^2  \phi  -\lambda  |\phi |^{p-1}\phi\,,
\end{equation}
where $ \mu >0$ and $\lambda >0 $. The equation for the Higgs real-valued scalar    field  in the 
  de~Sitter spacetime is a special case of (\ref{1.1}) when $p=3$, $n=3$:  
\begin{equation}
\label{Higgs_eq}
  \phi_{tt} +   3   \phi_t - e^{-2 t} \Delta  \phi =    \mu^2  \phi  -\lambda   \phi^ 3 \,.  
\end{equation}
Scalar fields play a fundamental role in the standard
model of particle physics, as well as its possible extensions.
In particular, scalar fields generate spontaneous
symmetry breaking and provide masses to gauge bosons
and chiral fermions by the Brout-Englert-Higgs mechanism~\cite{Englert-Brout}
using a Higgs-type potential~\cite{Higgs}.
\smallskip

The energy   
\begin{equation}
\label{finite_Energy}
E(t) =e^{nt}\int_{{\mathbb R}^n} \left( \frac{1}{2} |  \phi _t+ \frac{n}{2}\phi |^2+ \frac{1}{2}e^{-2 t}|\nabla_x   \phi|^2- \frac{1}{2}\left( \frac{n^2}{4} + \mu ^2 \right)|  \phi|^2 + \frac{1}{p+1}\lambda | \phi|^{p+1} \right) dx
\end{equation}
of the  $ L^q({\mathbb R}^n)  $-solution of the equation (\ref{1.1}) is non-increasing since
\[
\frac{d}{dt} E(t) =- e^{nt} \int_{{\mathbb R}^n} \left( e^{-2 t}|\nabla_x   \phi|^2 +  
 \lambda \frac{n(p-1)}{2(p+1)}  | \phi|^{p+1}\right) dx\,.
\]
The constants $\phi =  \pm \frac{\mu }{\sqrt{\lambda }} $ are non-trivial  real-valued solutions of the 
equation (\ref{Higgs_eq}) with the  positive  energy density   $e^{nt}\frac{\mu ^2}{8\lambda }(2n^2-1)$. The $x$-independent solution of  (\ref{Higgs_eq})
solves the Duffing's-type equation
\[
 \ddot \phi +3 \dot \phi =\mu^2  \phi  -\lambda   \phi^ 3\,,
\]
which describes the motion of a mechanical system in a twin-well potential field.
\smallskip

Unlike  the equation in the Minkowski spacetime, that is,  the equation 
\begin{equation}
\label{1.1_Min}
  \phi_{tt}  - \Delta  \phi =    \mu^2  \phi  -\lambda  |\phi |^{p-1}\phi\,,
\end{equation}
the equation (\ref{1.1}) has no other time-independent solution. For the  equation (\ref{1.1_Min}) the existence of a weak global solution in the energy space is known (see, e.g., 
 Proposition~3.2~\cite{G-V1989}) under certain conditions.
The equation 
\begin{equation}
\label{Higgs_eq_Min}
  \phi_{tt}      -   \Delta  \phi =    \mu^2  \phi  -\lambda   \phi^ 3 
\end{equation}
for the Higgs scalar    field  in the Minkowski spacetime has the time-independent flat solution 
\begin{equation}
\label{stadysol}
\phi_M (x) = \frac{\mu }{\sqrt{\lambda} } \tanh \left(   \frac{\mu^2 }{2 } N\cdot (x-x_0) \right), \qquad N, x_0, x \in {\mathbb R}^3.
\end{equation}
The unit vector $N $ defines the direction of  the propagation of the wave front. 
The  solution (\ref{stadysol}), after Lorentz transformation, gives rise
to  a traveling solitary wave of the form
\begin{equation}
\label{solitiry}
\phi_M (x,t) =  \frac{\mu }{\sqrt{\lambda} } \tanh \left(   \frac{\mu^2 }{2 }  [ N \cdot (x-x_0) \pm v(t-t_0)]\frac{1}{\sqrt{1-v^2}}\right), 
\qquad  N, x_0, x \in {\mathbb R}^3,
\,\, t \geq t_0,
\end{equation}
if $0<v<1$, where $v $ is the initial velocity. The set of zeros of the solitary wave $\phi =\phi_M (x,t) $, that is, the set given by $N \cdot (x-x_0) \pm v(t-t_0)=0 $, is the moving boundary of the {\it wall}.  Existence of standing waves $\phi =\exp(i\omega t)v(x) $,
which are exponentially small at infinity $|x| =\infty$, and of  corresponding  solitary waves  for the equation (\ref{Higgs_eq_Min}) 
with $  \mu^2 <0$ and $  \lambda  <0$ is known (see, e.g.,  \cite{Strauss}).
\smallskip

It is of  considerable interest for  particle physics and inflationary  cosmology  to study the so-called bubbles \cite{Coleman},
\cite{Linde}, \cite{Voronov}. 
In \cite{Lee-Wick} bubble is defined as a simply connected domain surrounded by a wall such that the field  
$ \phi (r,t)$ approaches one of the vacuums outside of a bubble. For a spherically symmetric bubble the  radius is defined at the zero of the function $\phi  $: 
$\phi (R(t),t)=0$. 
 The creation and
growth of bubbles is an interesting mathematical problem \cite[Ch.7]{Coleman},
\cite{Linde} and it motivates our interest in the sign changing  global solutions.
\smallskip 

A global  in time solvability of the Cauchy problem for equations (\ref{1.1}) and  (\ref{Higgs_eq})
 is not known. For the wave maps on Robertson--Walker   spacetimes, whose inverse radius is integrable with respect to the cosmic time,
 and, in particular, on the de~Sitter  spacetime,
a global existence is proven \cite{Choquet-Bruhat}.  
The local solution exists
 for every smooth initial data. 
The $C^2$ solution of the equation (\ref{Higgs_eq}) is unique and obeys the finite speed of the propagation 
property. (See, e.g., \cite{Hormander_1997}.) 
For the Higgs scalar    field  in the de~Sitter spacetime, that is for the equation (\ref{Higgs_eq}), Theorem~\ref{T_Higgs_boson} of Section~\ref{SHiggs} gives 
the necessary  conditions for the existence of the global solution. 
In particular, 
it gives the necessary conditions for the  global in time existence of bubbles.  In the forthcoming paper we prove the 
existence of the global in time small data solutions for the equation (\ref{1.1}). 
\smallskip

Since we are interested in the properties of global solutions
in the de~Sitter spacetime, we mention here two recent 
articles on linear equations on the asymptotically de~Sitter spacetimes. Vasy~\cite{Vasy_2010}
exhibited the well-posedness of the Cauchy problem and showed that on such spaces, the solution
of the Klein-Gordon equation without source term  and with smooth Cauchy data has an asymptotic expansion at infinity.  He also showed
that solutions of the wave equation exhibit scattering. Baskin~\cite{Baskin} constructed parametrix for the forward fundamental 
solution of the wave and Klein-Gordon equations on asymptotically de Sitter spaces without caustics and used this parametrix to obtain
asymptotic expansions for solutions of the equation with some class of source terms. 
(For more references on the asymptotically de Sitter spaces, see  the bibliography in \cite{Baskin},
 \cite{Vasy_2010}.) 
\smallskip

In order to make  our results   more transparent   
we formulate them for  the function \,$u = e^{\frac{n}{2}t}\phi$. For this  new unknown function \,$u=u(x,t) $,  the equation
(\ref{1.1})  takes the form of the semilinear 
Klein-Gordon  equation for \,$u$\, 
\begin{equation}
\label{2.2}
u_{tt} - e^{-2t} \bigtriangleup u  - M^2 u=  - \lambda e^{-\frac{n(p-1)}{2}t} |u |^{p-1}u,
\end{equation}
where the  ``curved mass'' \,$M\geq 0  $\, is defined as follows:
\begin{equation}
\label{M}
M^2:=  \frac{n^2}{4} + \mu ^2>  0\,.
\end{equation}
The equation (\ref{2.2}) is the so-called  equation with  imaginary mass. Equations with imaginary mass appear in 
several physical  models  such as \, $\phi ^4$ \, field model, tachion (super-light) fields,  Landau-Ginzburg equation  and others.
\smallskip

Next, we use the fundamental solution of the corresponding linear operator in order to reduce the Cauchy problem for the
semilinear equation to the integral equation and to define a weak solution. 
We denote by $G$ the resolving operator of the problem 
\begin{equation}
\label{2.2a}
u_{tt} - e^{-2t} \bigtriangleup u  - M^2 u=   f,   \quad u  (x,0) =  0 , \quad \partial_{t }u  (x,0 ) =0\,.
\end{equation}
Thus, $u=G[f]$. The  equation of (\ref{2.2a}) is strictly hyperbolic. This implies the well-posedness of the Cauchy problem (\ref{2.2a})  
in  different functional spaces. Consequently, the operator  $G$ is well-defined in those  functional spaces.
\smallskip

The operator $G$ is explicitly written   in \cite{Yag_Galst_CMP} for the case of the real mass; 
that is, with the mass term $+M^2u$ in the equation (\ref{2.2a}). 
The analytic continuation with respect to the parameter $M $  of this operator  allows us  also to use $G$
in the case of imaginary mass.  More precisely, for $M \geq 0$  
 we define the operator \,$G$\, acting on \,$f(x,t) \in C^\infty ({\mb R}\times [0,\infty)) $ \,by 
\begin{eqnarray*} 
G[f](x,t)   
& := &
 \int_{ 0}^{t} db \int_{ x - (e^{-b}- e^{-t})}^{x+ e^{-b}- e^{-t}}  dy\, f(y,b)  
(4e^{-b-t })^{-M} \Big((e^{-t }+e^{-b})^2 - (x - y)^2\Big)^{-\frac{1}{2}+M    } \\
&  &
\hspace{2.5cm}\times F\Big(\frac{1}{2}-M   ,\frac{1}{2}-M  ;1; 
\frac{ ( e^{-b}-e^{-t })^2 -(x- y )^2 }{( e^{-b}+e^{-t })^2 -(x- y )^2 } \Big)   , 
\end{eqnarray*}
where $F\big(a, b;c; \zeta \big) $ is the hypergeometric function. (See, e.g., \cite{B-E}.) For analytic continuation, see e.g., 
 \cite[Sec. 1.8]{Slater}. 
For $n \geq 2$, in both cases  of even and odd $n$, one can write 
\begin{eqnarray*} 
 G[f](x,t) 
&  =  &
2   \int_{ 0}^{t} db
  \int_{ 0}^{ e^{-b}- e^{-t}} dr  \,  v(x,r ;b)  (4e^{-b-t})^{-M}
\left( (e^{-t}  + e^{-b} )^2 - r^2   \right)^{-\frac{1}{2}+M}  \nonumber \\
&  &
\hspace{2.5cm} 
\times F\left(\frac{1}{2}-M,\frac{1}{2}-M;1; 
\frac{ (e^{-b}- e^{-t})^2-r^2}
{  (e^{-b}+ e^{-t})^2-r^2} \right) ,
\end{eqnarray*}
where the function 
$v(x,t;b)$   
is a solution to the Cauchy problem for the  wave equation
\[
v_{tt} -   \bigtriangleup v  =  0 \,, \quad v(x,0;b)=f(x,b)\,, \quad v_t(x,0;b)= 0\,.
\]
\medskip

It can be proved   that if $ n\big( \frac{1}{q'} - \frac{1}{q} \big) \leq 1$, $ \frac{1}{q'} +\frac{1}{q}=1$, $1\leq q' \leq 2 \leq  q \leq \infty $, then for every given $T>0$
the operator $G$ can be extended to a bounded operator:
\[
G \,: \, C([0,T]; L^{q'}({\mb R}^n)) \longrightarrow  C([0,T]; L^q({\mb R}^n)) \,.
\]
Consequently, the operator $G$ maps
\[
G \,: \, C([0,\infty); L^{q'}({\mb R}^n)) \longrightarrow  C([0,\infty); L^q({\mb R}^n)) ,
\]
in the corresponding topologies.  Moreover, 
\[
G \,: \, C([0,\infty); L^{q'}({\mb R}^n)) \longrightarrow  C^1([0,\infty); {\mathcal D}'({\mb R}^n))  .
\]

Let $u_0=u_0(x,t)$ be a solution of the Cauchy problem
\begin{equation}
\label{3.5}
\pa_{t }^2 u_0 - e^{-2t} \bigtriangleup u_0   - M^2 u_0 =   0, \quad u_0 (x,0) = \varphi _0(x), \quad \pa_{t }u_0 (x,0 ) = \varphi _1(x )\,.
\end{equation}
Then any solution $u=u(x,t)$ of the equation (\ref{2.2}), which takes initial value $ u (x,0) = \varphi _0(x), \quad \pa_{t }u (x,0) = \varphi _1(x)$,  
solves the integral equation
\begin{equation}
\label{2.6new}
u(x,t) = u_0(x,t)- G[\lambda e^{-\frac{n(p-1)}{2}\cdot }|u |^{p-1}u](x,t)   \,.
\end{equation}
We use the last equation   to  define  a weak solution of the problem for the differential equation.
Let $ \Gamma \in C([0,\infty))$. For every given function $u_0 \in C([0,T]; L^{q'}({\mb R}^n))$ we consider the integral equation  
\begin{equation}
\label{2.7}
u(x,t) = u_0(x,t)- G\left[\Gamma (\cdot ) \left| \int_{{\mb R}^n} |u (y,\cdot )|^{p-1} u (y,\cdot )\,dy \right|^\beta  |u (y,\cdot )|^{p-1} u (y,\cdot )\right](x,t)    \,,
\end{equation}
for the function 
\[
\dsp u  \in \bigcap_{i=1,p,q}C([0,T]; L^i({\mb R}^n)). 
\]
Here $\beta \in {\mathbb R} $, $ q' \geq  q>1$, $p  \geq 1$. The last integral equation 
corresponds to a slightly more general equation than (\ref{2.2}), namely, to the equation with {\it non-local nonlinearity} ({\it non-local self-interaction}) 
\begin{equation}
\label{2.8}
u_{tt} - e^{-2t} \bigtriangleup u  - M^2 u=  - \Gamma (t ) \left| \int_{{\mb R}^n} |u (y,t )|^{p-1}u (y,t ) dy \right|^\beta |u |^{p-1}u\,. 
\end{equation}
\begin{definition}
 If $u_0$ is a solution of the Cauchy problem (\ref{3.5}), then the solution  $u=u (x,t)$ of (\ref{2.7}) 
is said to be  
{\it a weak solution} of the Cauchy problem for the equation (\ref{2.8}) with the initial conditions
$ 
u (x,0) = \varphi _0 (x)$,\,   $\pa_{t }u (x,0) = \varphi _1 (x) .
$
\end{definition}

We are looking for  the  sufficient conditions for the continuous global solution for the changing of the sign, and, consequently,
for the creation of a bubble. It turns out that,  such conditions are also 
necessary conditions for the solution of equation (\ref{2.8}) to exist globally in time. 
 Equivalently,  
we are looking for  the  sufficient conditions on the  weak solution of the equation  that guarantee,   in general, the non-existence of
a  global in time weak solution, namely, the blow-up phenomena. 
In order to prove a sign-changing property of the  global solutions of the semilinear Klein-Gordon equations (\ref{1.1}), (\ref{1.1_Min})
 we invoke the so-called Functional Method that has been used to reveal blow-up phenomena for several equations (see, e.g., \cite[Ch. 2]{Alinhac}).
\smallskip

In the next definition we measure the variation of the sign of the function $\phi =\phi (x)  $ by the deviation from the H\"older inequality
of the inequality between the 
 integral of the function  and the self-interaction functional. 
Time $t$ is regarded as a parameter.
\begin{definition}
\label{D1.2}
The function $\phi  \in  C([0,\infty);$ $L^p(  {\mb R}^n ))$   
is said to be asymptotically time-weighted   $L^p$-non-positive  (non-negative), if  
there is a non-negative  number $C_\phi $ and positive non-decreasing 
function $\nu_\phi \in  C ([0,\infty)) $ such that with $\sigma =1$ ($\sigma =-1$)  one has 
\[
  \left| \int_{{\mathbb R}^n} \phi  (x,t) \, dx  \right|^{p }
\leq -\sigma C_\phi \nu_\phi  (t) \int_{{\mathbb R}^n}  |\phi  (x,t) |^{p-1}\phi (x,t)\, dx  
\quad \mbox{for all sufficiently large} \quad 
 t .
  \]
\end{definition}
It is evident  that any sign preserving function $\phi  \in   L^p(  {\mb R}^n ) $ with a compact support 
satisfies the last inequality with $ \nu_\phi  (t) \equiv 1 $  
and either $\sigma =1 $ or $\sigma =-1$, while $C_\phi^{ 1/( p-1) } $ is a measure of the support.
Then, any smooth global non-positive (non-negative)  solution $\phi = \phi (x,t) $ of the Cauchy problem for (\ref{1.1}) or (\ref{1.1_Min}) with compactly supported initial data is
also asymptotically time-weighted   $L^p$-non-positive (non-negative) with the weight $\nu_\phi  (t)= (1+t)^{n(p-1)} $.

\begin{theorem}
\label{T2.1}
Let $u=u(x,t) \in C([0,\infty);L^q(  {\mb R}^n ))$, $2\leq q < \infty $, be a global solution
of the equation (\ref{2.7}) with $ \beta >1/p-1 $. Suppose that the function  $ \Gamma \in C^1([0,\infty))$ is either non-decreasing or non-increasing. 
Assume that   
 the function $u_0 \in C^1({\mb R}^n \times [0,\infty))$ satisfies 
\begin{eqnarray}
\hspace{-1.3cm} &  &
\label{21}
\sigma \left( M\int_{{\mathbb R}^n}  u_0(x,0)  \, dx + \int_{{\mathbb R}^n}  \partial_t u_0(x,0)  \, dx \right)>0\,,
\end{eqnarray}
and, additionally,
\begin{eqnarray}
\hspace{-1.3cm} &  &
\label{2.9}
\int_{{\mathbb R}^n} u_0(x,t)  dx =   \cosh (Mt)\int_{{\mathbb R}^n}  u_0(x,0)  \, dx 
+ \frac{1}{M} \sinh(M t) \int_{{\mathbb R}^n}  \partial_t u_0(x,0)  \, dx
\quad \mbox{for all} \quad t \geq 0
\end{eqnarray}
if $M>0$, while for $M=0$
\begin{eqnarray}
\hspace{-1.3cm} &  &
\label{2.9_0}
\int_{{\mathbb R}^n} u_0(x,t)  dx =   \int_{{\mathbb R}^n}  u_0(x,0)  \, dx + t \int_{{\mathbb R}^n}  \partial_t u_0(x,0)  \, dx
\quad \mbox{for all} \quad t \geq 0.
\end{eqnarray}
Assume also that the self-interaction functional satisfies
\begin{eqnarray}
\label{17}
 &  &
\sigma  \int_{{\mathbb R}^n}   |u (z,b) |^{p-1}u(z,b)\, dz  
\leq 0   
\end{eqnarray}
for all \, $t$\,  either outside of the sufficiently small neighborhood of zero 
if \,$M  >0$, or  inside  of some neighborhood of infinity 
if \,$M =0$.
 
Then, the global solution $u=u(x,t)$ cannot be 
an  asymptotically   time-weighted $L^p$-non-positive (-non-negative)  with the weight  $\nu_u \in  C^1([0,\infty)) $ such that 
 if $M> 0$, then 
\beqst
&  &
\Gamma (t) \geq c \nu_u  (t)^{\beta +1}e^{-M(p(\beta +1)-1 )t} t^{2+\varepsilon}\quad \mbox{for all large} \quad t    
\eeqst
with the numbers $\varepsilon >0$ and $ c>0$, while for $M=0$ it satisfies
\beqst
&  &
\Gamma (t) \geq c t^{-1-p(\beta +1)} \nu _u(t)^{\beta +1}\quad \mbox{for all large} \quad t\,.
\eeqst
\end{theorem}

 The conditions (\ref{2.9}), (\ref{2.9_0}) are inherited from the partial differential equation (\ref{3.5}). 
In fact, every smooth integrable solution of (\ref{3.5}) satisfies (\ref{2.9}) or (\ref{2.9_0}). 

Thus, the theorem shows that the continuous,  asymptotically   time-weighted  $L^p$-non-positive  
(non-negative)  global solution of the equation  (\ref{2.7}) cannot be sign preserving  
if it is generated by the function $u_0=u_0(x,t) $, which obeys  (\ref{21}), (\ref{2.9}), (\ref{2.9_0}). 
Consequently, smooth  asymptotically   time-weighted  $L^p$-non-positive  
(non-negative)  global solution of the equation  (\ref{2.2}) cannot be sign preserving  
if its initial data $\varphi _0 $, $\varphi _1 $ imply  (\ref{21}).

An application of the last theorem to the 
generalized Higgs real-valued scalar field 
equation (\ref{1.1}) 
with $\mu >0$ results in the following corollary.
\smallskip

 \begin{corollary}
Let $\phi =\phi (x,t) \in C([0,\infty);L^q(  {\mb R}^n ))$, $2\leq q < \infty $, be a global   weak solution
of the equation  (\ref{1.1}).
Assume also that
 the initial data of $\phi =\phi (x,t)$ satisfy 
\begin{eqnarray*}
\sigma \left( \left(\sqrt{\frac{n^2}{4} + \mu ^2} +\frac{n}{2}\right)C_0(\phi )+C_1(\phi ) \right) >0  
\end{eqnarray*}
with $\sigma =1$ ($\sigma =-1$), while 
\begin{eqnarray*} 
\sigma  \int_{{\mathbb R}^3}   |\phi (x,t)|^{p-1}\phi (x,t)\, dx  \leq  0 
\end{eqnarray*}  
is fulfilled 
 for all\, $t$\,  outside of the sufficiently small neighborhood of zero.

Then, the global solution $\phi =\phi (x,t)$ cannot be an  asymptotically   time-weighted $L^p$-non-positive (-non-negative)  solution  with the weight  
$\nu_\phi   (t)=e^{a_\phi  t }t^{b_\phi } $, where either 
$ a_\phi <(p-1)\left( \sqrt{\frac{n^2}{4}+\mu ^2}- \frac{n }{2}\right)$, $b_\phi \in{\mathbb R} $, or 
$   a_\phi = (p-1)\left( \sqrt{\frac{n^2}{4}+\mu ^2}- \frac{n }{2}\right)$, $b_\phi <-2 $. 
 \end{corollary}

The numbers $\pm \sqrt{\frac{n^2}{4}+\mu ^2}-\frac{n }{2}$ of the last corollary are the roots  of the characteristic equation of the linear ordinary differential part of 
(\ref{1.1}). They  also appear in the 
 exponents $s_\pm(\lambda )$  of  the suggested in   \cite{Vasy_2010}  representation of the solution of the 
Cauchy problem for the linear Klein-Gordon operator on the asymptotically de~Sitter-like spaces.

 \smallskip

Consider now the case of the  function $u_0=u_0(x,t) $, which has support in 
the cylinder $B^n_R(0)\times [0,\infty)$, where $B^n_R(0)$ is a ball in ${\mathbb R}^n $ with radius $R>0$.
The smooth solution of the Cauchy problem for the linear equation (\ref{3.5}) 
with the compactly supported initial data can exemplify such $u_0=u_0(x,t) $.
The diversity of the eigenfunctions of the Laplace operator  gives a more sensitive test 
to find out the sign changing solutions for the  equation (\ref{2.8}) with $\beta =0$.
 \begin{definition}
 \label{D3}
The function $\phi  \in  C([0,\infty);L^p(  {\mb R}^n ))$ 
is said to be asymptotically time-weighted  $-\psi$~$L^p$ -signed   if  
there are  
 a non-negative  number $C_{\phi, \psi } $, and positive 
non-decreasing function $\nu_{\phi, \psi } \in  C ([0,\infty)) $ 
such that   the following inequality holds
\[
  \left| \int_{{\mathbb R}^n}  \psi (x) \phi  (x,t) \, dx  \right|^{p }
\leq - C_{\phi, \psi } \nu_{\phi, \psi }  (t) \int_{{\mathbb R}^n}  \psi (x)  
|\phi  (x,t) |^{p-1}\phi (x,t)\, dx  \quad \mbox{for all large} \quad t.
  \]
\end{definition}
It is evident  that if with some  eigenfunction $\psi  $  the  function $\psi \phi  \in   L^p(  {\mb R}^n ) $ 
has a compact support and is  sign preserving, then 
 the last inequality holds with  either $ \nu_{\phi, \psi }  (t) \equiv 1 $  
 or $ \nu_{\phi, \psi }  (t) \equiv -1 $, while $C_{\phi, \psi }^{1/( p-1)} $ is a measure of the support.

  \begin{theorem}
\label{T2.1_psi}
Let $u=u(x,t) \in C([0,\infty);L^q(  {\mb R}^n ))$, $2\leq q < \infty $, be a global solution
of the equation 
\begin{equation}
\label{2.7b0}
u(x,t) = u_0(x,t)- G\left[\Gamma (\cdot )  |u (y,\cdot )|^{p-1} u (y,\cdot )\right](x,t)    
\end{equation}
with  $p >1$. Suppose that the function  $ \Gamma \in C^1([0,\infty))$ is either non-decreasing or non-increasing. 
  Let $\psi =\psi (x)$ be an eigenfunction  $\psi =\psi (x) $ of the Laplace operator in ${\mathbb R}^n $
corresponding to the eigenvalue $\nu  $. Further, suppose that 
$u_0 \in C^1({\mb R}^n \times [0,\infty))$ with the support in the $B^n_R(0)\times [0,\infty)$ satisfies 
\beq
 &  &
\label{21_psi}
 (M +\nu )\int_{{\mathbb R}^n} \psi (x) u_0(x,0)  \, dx + \int_{{\mathbb R}^n}   \psi (x) \partial_t u_0(x,0)  \, dx >0\,, 
\eeq
and, additionally
\beq
 &  &
\label{2.9_psi}
\int_{{\mathbb R}^n} \psi (x)  u_0(x,t)  dx =   \cosh ((M +\nu )t)\int_{{\mathbb R}^n}  \psi (x)  u_0(x,0)  \, dx \\
&  &
\hspace{3.3cm} + \frac{1}{M+\nu } \sinh((M +\nu ) t) \int_{{\mathbb R}^n}   \psi (x) \partial_t u_0(x,0)  \, dx
\quad \mbox{for all}\quad t \geq 0 , \nonumber
\eeq
if $M +\nu >0 $, while for $M +\nu =0 $
\beq
 &  &
\label{2.9_psi_0}
\int_{{\mathbb R}^n} \psi (x)  u_0(x,t)  dx =    \int_{{\mathbb R}^n}  \psi (x)  u_0(x,0)  \, dx  +   t \int_{{\mathbb R}^n}   \psi (x) \partial_t u_0(x,0)  \, dx
\quad \mbox{for all}\quad t \geq 0 \,. 
\eeq
  Assume that the self-interaction functional satisfies
\begin{eqnarray*}
 &  &
\int_{{\mathbb R}^n}  \psi (x)  |u (z,b) |^{p-1}u(z,b)\, dz 
   \leq 0   \,,
\end{eqnarray*}
for all $t$ either outside of the sufficiently small neighborhood of zero 
if \,$M +\nu  >0$, or  inside  of some neighborhood of infinity
if \,$M +\nu =0$.

Then, the global solution $u=u(x,t)$ cannot be 
an  asymptotically   time-weighted   $-\psi$~$L^p$-signed  with the weight  $\nu_{u,\psi } \in  C^1([0,\infty)) $ such that 
  if $M+\nu > 0$, then 
\beqst
&  &
 \Gamma (t) \geq c \nu_{u,\psi }  (t)e^{-(M +\nu )(p-1 )t} t^{2+\varepsilon}\quad 
\mbox{for all large} \quad t    
\eeqst
with the numbers $\varepsilon >0$ and $ c>0$, while for $M+\nu =0$ it satisfies
\beqst
&  &
 \Gamma (t) \geq c t^{-1-p} \nu _{u,\psi }(t)\quad \mbox{for all large} \quad t \,.
\eeqst 

\end{theorem}
\smallskip

For the differential equation (\ref{Higgs_eq}) of the  Higgs boson in the de~Sitter spacetime, Theorem~\ref{T2.1_psi} leads to the following result   
 (see Teeorem~\ref{T_Higgs_boson}):  
 The continuous global solutions  obtained by prolongation of some local solutions must change a sign, 
and consequently, they vanish    at some points. In particular, such radial global solutions have zeros and, therefore, they
give rise
to  at least one bubble.   
Hence, for the global solutions Theorem~\ref{T_Higgs_boson} guarantees the  creation of the bubble. 
Moreover, acccording to Corollary~\ref{C4.4}
 the bubbles exist in any neighborhood of infinite time. 
 Thus, the global solution is an oscillating in time solution. 
In particular, for the  continuous global solutions we give integral  conditions (see below (\ref{sign})  and (\ref{C0C1})), which are 
 sufficient conditions for the creation of the bubbles and their existence in the future. 
Similar conclusions  are valid for the equations in the Minkowski spacetime (Section~\ref{S_Mink}).
 \smallskip

\medskip

This paper is organized as follows.  In Section \ref{S_Mink} we discuss and prove properties of the  Higgs boson equation in the 
Minkowski spacetime.
First in Theorem~\ref{T_Higgs_boson_Mink} we give some criteria via integrals of the solution, which has   unbounded support in the spatial variables.  
Next we use the eigenfunctions of the Laplace operator to widen those criteria  by weighted integrals of the solutions, 
which has  compact supports in the spatial variables. 
In Section \ref{S_Proofs} we prove Theorem~\ref{T2.1} about solutions of the  
equation in the de~Sitter spacetime with, in general,   non-local self-interaction.  In Section \ref{SHiggs}
we prove  Theorem~\ref{T2.1_psi} and then in Theorems~\ref{T_Higgs_boson}-\ref{T_Higgs_boson_Lap} we discuss  in more detail  the case of solutions of the Higgs boson equations 
in the de~Sitter spacetime with $x \in {\mathbb R}^3$, which have  compact supports in the spatial variables.   
In the Appendix, Section~\ref{Appendix}, we give some integral representations of the hyperbolic sine function and one generalization of 
Kato's lemma to the second order ordinary differential inequality with the exponentially decaying  kernel.

\medskip

\section{The Higgs boson in the Minkowski spacetime} 
\label{S_Mink}

Consider the generalized Higgs boson  equation in the Minkowski spacetime, 
\begin{equation}
\label{1.1_Min_Mink}
  \phi_{tt}  - \Delta  \phi =    \mu^2  \phi  -\lambda  |\phi |^{p-1}\phi\,.
\end{equation}
Here for the numbers $\mu  $ and $\lambda  $ we assume $\mu \geq  0 $ and $\lambda >0 $.
For the  equation (\ref{1.1_Min_Mink}) the existence and the uniqueness of weak global solution  $(\phi ,\dot \phi  ) \in C({\mathbb R};X_e)$
in the energy space 
$ X_e:=H^1\oplus  L^2$ is known (see,  
 Proposition~3.2~\cite{G-V1989}  and Theorems~6.2-6.3~\cite{Shatah}) under certain conditions on $n$ and $p$, which 
include $p-1<4/(n-2) $. That solution satisfies the conservation of energy $
E(\phi (t),\dot \phi (t))= E(\phi  (0),\dot \phi  (0)) \equiv E$
 and the following estimates
 \begin{eqnarray*}
  &  &
\int_{{\mathbb R}^{n}} |\phi (x,t)|^2\,dx \leq e(E,t)\,,\\
&  &
\int_{{\mathbb R}^{n}} \left( |\dot \phi (x,t)|^2 + |\nabla \phi (x,t)|^2 \right)\,dx 
\leq \dot e(E,t)\,,
\end{eqnarray*} 
where
\[
   e(E,t):=  \left(\int_{{\mathbb R}^{n}} |\phi (x,0)|^2\,dx \right)\cosh (\mu t )
  + \left(E+\mu ^2\int_{{\mathbb R}^{n}} |\phi (x,0)|^2\,dx\right)^{1/2} \mu ^{-1} \sinh (\mu t )\,.
 \]
Moreover, if $p<\frac{n+2}{n-2}$ and $n \leq 9$, then for smooth data 
$\phi  (0), \dot \phi  (0) \in C^\infty$, the solution $\phi $ is $C^\infty$-smooth.
The equation 
\begin{equation}
\label{Higgs_eq_Min_Mink}
  \phi_{tt}      -   \Delta  \phi =    \mu^2  \phi  -\lambda   \phi^ 3 
\end{equation}
for the Higgs scalar    field  in the Minkowski spacetime has the  time-independent real-valued flat solution (\ref{stadysol}) as well as a traveling solitary wave.
\smallskip 

In  this paper we are looking for  the qualitative properties of the global in time solutions of equation (\ref{1.1_Min_Mink}).
More precisely, we are interested in the sign changing   global solutions of equations (\ref{1.1_Min_Mink}) and (\ref{Higgs_eq_Min_Mink}).
Our interest in   such solutions   is motivated by the problem of  the creation and
growth of bubbles.
The   theorems  below give 
 necessary conditions for the   global  in time existence of the bubbles. 
\smallskip 

In order to solve the Cauchy problem for the semilinear equation via the integral 
equation and to define a weak solution, we use the fundamental solution of the corresponding linear operator. 
We denote by $G$ the resolving operator of the problem 
\begin{equation}
\label{2.2a_Mink}
u _{tt} -  \bigtriangleup u   - \mu ^2 u =   f,   \quad u   (x,0) =  0 , \quad \partial_{t }u   (x,0 ) =0\,.
\end{equation}
Thus, $u=G[f]$. One can write the following explicit formula (see, e.g., \cite{yagdjian_Rend_Trieste}) for the operator $G$, namely,
\begin{eqnarray}
\label{0.8K-G-Im}
\hspace{-0.5cm} u(x,t)
\!\!& \!\!= \!\!&\!\!
\int_{ 0}^{t} db
  \int_{ 0}^{t-b }   I_0\left(\mu \sqrt{(t-b)^2-r^2} \right)  w(x,r;b )  dr
, \,\, x \in {\mathbb R}^n,
\end{eqnarray}
where $w(x,t;b ) $ is a solution of  
\[
\cases{   w_{tt}-\Delta w = 0\quad  \mbox{\rm in} \quad  {\mathbb R}^{n+1}, \cr
w(x,0 ;\tau )=f(x,\tau ),\quad  w_t(x,0 ;\tau )= 0  \quad  \mbox{\rm in} \quad  {\mathbb R}^{n},  }
\]
while $\tau $ is regarded as a parameter. The function $I_0(z) $ is the modified Bessel function of the first kind.
For $w= w(x,t;b) $ there are the following  representation formulas (see, e.g., \cite{Shatah}).  If $n$ is odd, $n=2m+1$,  $m  \in {\mathbb N}$,  then for $f\in C^\infty  ({\mb R}^n\times [0,\infty)) $, we have
\begin{eqnarray*} 
w(x,t;b)  
& = &
\frac{\partial }{\partial r} 
\Big(  \frac{1}{r} \frac{\partial }{\partial r}\Big)^{\frac{n-3}{2} } 
\frac{r^{n-2}}{\omega_{n-1} c_0^{(n)} }  \!\!\int_{S^{n-1} } f(x+ry,b)\, dS_y    ,
\end{eqnarray*} 
where $c_0^{(n)} =1\cdot 3\cdot \ldots \cdot (n-2 )$.
The constant $\omega_{n-1} $ is the area of the unit sphere $S^{n-1} \subset {\mathbb R}^n$. 

If $n$ is even, $n=2m$,  $m  \in {\mathbb N}$,  then for $f  \in C^\infty   ({\mb R}^n\times [0,\infty)) $, we have
\begin{eqnarray*}
w(x,t;b)   
& = &
\frac{\partial }{\partial r} 
\Big( \frac{1}{r} \frac{\partial }{\partial r}\Big)^{\frac{n-2}{2} } 
\frac{2r^{n-1}}{\omega_{n-1} c_0^{(n)} }  \!\!\int_{B_1^{n}(0)} \frac{f(x+ry,b) }{\sqrt{1-|y|^2}} \, dV_y .
\end{eqnarray*} 
Here $B_1^{n}(0) :=\{|y|\leq 1\} $ is the unit ball in ${\mathbb R}^n$, while $c_0^{(n)} =1\cdot 3\cdot \ldots \cdot (n-1)$.
\medskip

In particular,
\begin{equation}
\label{conslaw}
\int_{{\mathbb R}^n}  w(x,t;b )\,dx= \int_{{\mathbb R}^n}  f(x, b )\,dx \quad \mbox{\rm for all}\quad  t\,.
\end{equation}
 The  equation of (\ref{2.2a_Mink}) is strictly hyperbolic and the operator  $G$ is well-defined in the several  functional spaces.

Let $\phi _0=\phi _0(x,t)$ be a solution of the Cauchy problem
\begin{equation}
\label{3.5_Mink}
\pa_{t }^2 \phi _0 -   \bigtriangleup \phi _0   - \mu ^2 \phi _0 =   0, \quad \phi _0 (x,0) = \varphi _0(x), \quad \pa_{t }\phi _0 (x,0 ) = \varphi _1(x )\,.
\end{equation}
Then any solution $\phi =\phi (x,t)$ of the equation (\ref{1.1_Min_Mink}), which takes initial value $ \phi  (x,0) = \varphi _0(x), \quad \pa_{t }\phi  (x,0) = \varphi _1(x)$,  
solves the integral equation
\begin{equation}
\label{2.6new_Mink}
\phi (x,t) = \phi _0(x,t)- G[\lambda |\phi  |^{p-1}\phi ](x,t)   \,.
\end{equation}
For every given function $\phi _0 \in C([0,T]; L^{q'}({\mb R}^n))$ we consider the integral equation  (\ref{2.6new_Mink}) 
for the function 
\[
\dsp \phi   \in \bigcap_{i=1,p,q}C([0,T]; L^i({\mb R}^n)). 
\]
Here $ q' \geq  q>1$, $p  \geq 1$. 
 If $\phi _0$ is generated by the Cauchy problem (\ref{3.5_Mink}), then the solution  $\phi =\phi  (x,t)$ of (\ref{2.6new_Mink}) 
is said to be  
{\it a weak solution} of the Cauchy problem for equation (\ref{1.1_Min_Mink}) with the initial conditions
$ 
\phi  (x,0) = \varphi _0 (x)$,\,   $\pa_{t }\phi  (x,0) = \varphi _1 (x) .
$

\begin{theorem}
\label{T_Higgs_boson_Mink}
Let $\phi=  \phi (x,t)\in C ({\mathbb R}^n\times [0,\infty))$ be a weak global solution   of  the
real  field  equation (\ref{1.1_Min_Mink}). Denote the integrals (functionals) of the initial values  
of $\phi $ by
\begin{equation}
\label{2.10_Mink}
C_{0   } (\phi):= \int_{{\mathbb R}^n}  \phi  (x,0 )   dx,  \quad
C_{1   }(\phi):= \int_{{\mathbb R}^n }  \phi_t (x,0 )   dx \,.
\end{equation} 
Assume that
the   self-interaction functional \,\,
 $-\lambda\int_{{\mathbb R}^n}  |\phi (x,t)|^{p-1}\phi (x,t)\, dx$   satisfies
 \begin{eqnarray}
\label{sign_Mink}
 \sigma  \int_{{\mathbb R}^n}  |\phi (x,t)|^{p-1}\phi (x,t)\, dx  \leq  0   \,,
\end{eqnarray} 
for all $t$ either outside of the sufficiently small neighborhood of zero 
if \, $\mu  >0$, or  inside  of some neighborhood of infinity
if \,$\mu =0$.  Assume also  that 
\begin{eqnarray}
\label{C0C1_Mink}
\sigma  \left( \mu C_{0   }(\phi )+C_{1   }(\phi ) \right) >0\,, 
\end{eqnarray}
where  either $\sigma  =-1$ or $\sigma  =1$. 

Then, the global solution $\phi =\phi (x,t)$ cannot be 
an  asymptotically   time-weighted $L^p$-non-positive (-non-negative)  with the weight  $\nu _\phi = e^{a_\phi t} t^{b_\phi }$, where  
 if  $\mu >0$, then  either  $  a_\phi < \mu (p-1)$  or $ a_\phi = \mu (p-1)$ and 
$b_\phi <-2$, while $a_\phi =0 $ and $b_\phi \leq 1+p $ if $\mu = 0$. 
\end{theorem}
\medskip

\noindent
{\bf Proof.} 
We consider the case of $\sigma =1$ only, since the case of  $\sigma =-1$ follows by the reflection $\phi \rightarrow  -\phi  $. 
We discuss separately   two  cases:  with  positive  mass, $\mu >0$, and vanishing  mass, $\mu =0$, respectively.
We start with the case of positive  mass. Let 
$\phi _0 \in C^1([0,\infty)\times {\mb R}^n )$ be a function with 
\begin{equation} 
\label{u0cond}
\int_{{\mathbb R}^n} \phi _0(x,t)  dx = C_0 (\phi)\cosh (\mu t) +C_1(\phi)\frac{1}{\mu } \sinh(\mu  t) 
\qquad \mbox{\rm for all} \quad t \geq 0\, , 
\end{equation}
where
\begin{equation}
\label{38}
\phi _0(x,0 )= \phi  (x,0 ),  \quad \partial_t \phi _0(x,0 ) = \partial_t \phi  (x,0 ) \,.
\end{equation}
The integrable in ${\mathbb R}^n_x$ solution to the problem (\ref{3.5_Mink}) satisfies (\ref{u0cond}). 
Thus, $\phi   \in C([0,\infty); L^q({\mb R}^n))$  is a solution to 
(\ref{2.6new_Mink}) generated by $\phi _0 $. 
According to the definition of the solution, for every given $T>0$ we have 
\[
G\left[  |\phi  |^{p-1}\phi \right]  \in C ([0,T]; L^{q }({\mathbb R}^n))\cap C ^1([0,T]; {\mathcal D}'({\mb R}^n))
\]
and (\ref{38}).
Then  $\phi $ is a continuous function of $t \in [0,\infty)$ with values in  $\ L^1({\mb R}^n)$, and 
we may integrate the equation (\ref{2.6new_Mink}):
\begin{equation}
\label{1.7a_Mink} 
\int_{{\mb R}^n} \phi (x,t)  \,dx = \int_{{\mb R}^n} \phi _0(x,t)  \,dx 
- \lambda \int_{{\mb R}^n } G\left[  |\phi  |^{p-1}\phi \right] (x,t)  \,dx   .
\end{equation}
In particular,
\begin{equation}
\label{40} 
 \int_{{\mb R}^n} \phi  (x,0) \,dx= \int_{{\mathbb R}^n} \phi _0(x,0 ) \,  dx=C_0, 
\quad \int_{{\mb R}^n} \phi _t (x,0)\,dx = \int_{{\mathbb R}^n } \partial_t \phi _0(x,0 ) \,  dx=C_1\,.
\end{equation}
Then, 
for the smooth function $\phi =\phi (x,t) $     we obtain 
\begin{eqnarray*}
\hspace{-0.5cm}
\int_{{\mathbb R}^n} G\left[  |\phi  |^{p-1}\phi \right](x,t) \, dx  
 &  =  & 
\int_{{\mathbb R}^n} \, dx\int_{ 0}^{t} db
  \int_{ 0}^{t-b }   I_0\left(\mu \sqrt{(t-b)^2-r^2} \right)  w(x,r;b )  dr\,,
\end{eqnarray*}
 where $w(x,r;b ) $ is solution of the problem
\[
\cases{   w_{tt}-\Delta w = 0\quad  \mbox{\rm in} \quad  {\mathbb R}^{n+1}, \cr
w(x,0 ;\tau )=|\phi (x,\tau ) |^{p-1}\phi (x,\tau ),\quad  w_t(x,0 ;\tau )= 0  \quad  \mbox{\rm in} \quad  {\mathbb R}^{n}\, , }
\]
with the parameter $\tau \in [0,\infty) $. Therefore, in view of (\ref{conslaw}),
\begin{eqnarray*}
\hspace{-0.5cm}
\int_{{\mathbb R}^n} G\left[  |\phi  |^{p-1}\phi \right](x,t) \, dx  
 &  =  & 
\int_{ 0}^{t} db
\int_{ 0}^{t-b }   I_0\left(\mu \sqrt{(t-b)^2-r^2} \right)  dr 
\int_{{\mathbb R}^n} dx\, w(x,0;b ) \\
& = &
\int_{ 0}^{t} db
  \int_{ 0}^{t-b }   I_0\left(\mu \sqrt{(t-b)^2-r^2} \right)  dr \int_{{\mathbb R}^n}   |\phi (x, b ) |^{p-1}\phi (x, b )\,dx  \\
 & = &
\int_{ 0}^{t} db \int_{{\mathbb R}^n}   |\phi (x, b ) |^{p-1}\phi (x, b )\,dx 
  \int_{ 0}^{t-b }   I_0\left(\mu \sqrt{(t-b)^2-r^2} \right)  dr  \,.
\end{eqnarray*}
One can easily check that 
\begin{eqnarray*}
  \int_{ 0}^{t-b }   I_0\left(\mu \sqrt{(t-b)^2-r^2} \right)  dr  =\frac{1}{\mu }\sinh(\mu (t-b))\,.
\end{eqnarray*}
Thus, we obtain 
 \begin{eqnarray*} 
\int_{{\mathbb R}^n} G[|\phi  |^{p-1}\phi  ](x,t)  dx  
  & = &
\int_{ 0}^{t}     \left( \int_{{\mathbb R}^n}  |u (x,b) |^{p-1}u(x,b)\, dx \right) 
 \frac{1}{\mu } \sinh (\mu (t-b ) ) \,  db .
\end{eqnarray*} 
Hence, (\ref{1.7a_Mink}) reads as follows:
\begin{eqnarray*} 
\int_{{\mathbb R}^n} \phi (x,t)  \,dx 
& = &
\int_{{\mathbb R}^n} \phi _0(x,t)  \,dx 
- \lambda  \int_{ 0}^{t}   \left( \int_{{\mathbb R}^n}  |\phi  (x,b) |^{p-1}\phi (x,b)\, dx \right)
 \frac{1}{\mu } \sinh (\mu (t-b ) ) \,  db \, .  
\end{eqnarray*}
Taking into account (\ref{40})   we derive
\begin{eqnarray*}
\int_{{\mathbb R}^n} \phi (x,t)  \,dx 
& = &
\frac{1}{2}\left(C_0+\frac{C_1}{\mu } \right)e^{\mu  t} +\frac{1}{2}\left(C_0-\frac{C_1}{\mu } \right) e^{-\mu  t}  \nonumber\\ 
&  &
-  \lambda \int_{ 0}^{t}  \left( \int_{{\mathbb R}^n}  |\phi  (x,b) |^{p-1}\phi (x,b)\, dx \right) 
 \frac{1}{\mu } \sinh (\mu (t-b ) ) \,  db   \, .
\end{eqnarray*}
Thus,   we obtain 
\begin{eqnarray*}
F(t)
& = &
   C_0 \cosh(\mu t)+ \frac{C_1}{\mu }  \sinh (\mu t) \\
&  &  
-  \lambda \int_{ 0}^{t}  \left( \int_{{\mathbb R}^n}  |\phi (x,b) |^{p-1}\phi (x,b)\, dx \right) 
 \frac{1}{\mu } \sinh (\mu (t-b ) ) \,  db    \,,
\end{eqnarray*}
where we have denoted 
\begin{eqnarray}
\label{F}
F(t)
& := &
  \int_{{\mathbb R}^n } \phi (x,t)   \,dx     \,.
\end{eqnarray}
It follows $ F \in C^2([0,\infty))$. 
Moreover, 
\begin{eqnarray}
\label{1.17a_Mink}
\dot F(t) 
& = &
   C_1\cosh(\mu t)+ \mu C_0  \sinh (\mu t)   \\
&  &   
-  \lambda \int_{ 0}^{t} \left( \int_{{\mathbb R}^n}  |\phi (x,b) |^{p-1}\phi (x,b)\, dx \right)
\cosh (\mu (t-b ) ) \,  db   \,, \nonumber \\ 
\ddot F(t) 
& = &
 \mu ^2F(t)   -     \lambda  \int_{{\mathbb R}^n}  |\phi (x,t) |^{p-1}\phi (x,t)\, dx    \,.
\end{eqnarray}
In particular, since (\ref{sign_Mink}) and  $ C_\phi  \geq 0$, $\nu_\phi   (t) \geq 0 $, there is a 
positive number $\varepsilon  $, $0<\varepsilon <1$,  such that
\begin{eqnarray}
\label{3.21new_Mink}
F(t)
& \geq &
 (1-\varepsilon)\left(   C _0    \cosh(\mu t)+ \frac{C_1}{\mu }  \sinh (\mu t) \right)\qquad 
\mbox{\rm for large} \,\, \, t \,.   
\end{eqnarray}
Indeed, according to (\ref{sign_Mink}) there exist positive $\varepsilon <1$ and $\delta_\phi>0 $ such that
\begin{equation}
\label{nbdzero}
 \lambda  \left|  \int_{ 0}^{\delta_\phi  }  \left( \int_{{\mathbb R}^n}  |\phi (x,b) |^{p-1}\phi (x,b)\, dx \right) 
 \frac{1}{\mu } \sinh (\mu (t-b ) ) \,  db   \right| \leq \varepsilon \cosh(\mu t) \quad 
\mbox{\rm for large} \,\, \, t \,.
\end{equation}
Then, the inequality (\ref{3.21new_Mink}) is fulfilled for  all    $ t \geq \delta_\phi $, if $\delta_\phi $ is sufficiently large.
By means of the   condition $\mu  C _0  + C_1   >0 $    we   conclude that 
\[
F(t) \geq 0  \qquad 
\mbox{\rm for large} \,\, \, t  \,,
\]
and, consequently,
\begin{eqnarray}
\ddot F(t)  
& \geq &
-  \lambda   \int_{{\mathbb R}^n}  |\phi  (x,t) |^{p-1}\phi (x,b)\, dx   \qquad 
\mbox{\rm for large} \,\, \, t  .
\end{eqnarray}
On the other hand, using Definition~\ref{D1.2}  with $\nu_\phi   (t) =e^{a_\phi t} t^{b_\phi } $   we obtain
\begin{eqnarray*}
\left| \int_{{\mathbb R}^n} \phi  (x,t) \,dx \right|^{p}
& \le &
-  C_\phi  \nu_\phi   (t)
\int_{{\mathbb R}^n}  |\phi  (x,t) |^{p-1}\phi (x,t)\, dx \\
& \leq  &
\lambda ^{-1}C_\phi  \nu_\phi   (t)  
 \left(- \lambda    \int_{{\mathbb R}^n} |\phi (x,t)|^{p-1}\phi (x,t) \,dx \right) \\
& = &
\lambda ^{-1} C_\phi   \nu_\phi   (t)    \Big(\ddot{F}(t)-  \mu ^2F(t) \Big)  \\
& \leq  &
\lambda ^{-1} C_\phi  \nu_\phi   (t)   \ddot{F}(t)   \qquad 
\mbox{\rm for large} \,\, \, t  \,.
\end{eqnarray*}
Thus, since $\nu_\phi   (t)>0 $ we have 
\[
\ddot{F}(t) 
\ge  \delta _0 \nu_u  (t)^{ -1}  F(t) ^{p}\qquad 
\mbox{\rm for all large} \,\, \, t \quad \mbox{\rm with} \,\, \,\delta _0 :=\lambda C_\phi ^{-1} >0 \,.
\]
Hence,    the last inequality together with (\ref{1.17a_Mink}) to (\ref{3.21new_Mink})
 implies the following system of the ordinary differential inequalities 
\[ 
  \left\{ 
\begin{array}{ccccc}
\dsp F(t) & \geq  & \hspace{- 0.2cm}  (1-\varepsilon ) C_0 \cosh(\mu t)+ (1-\varepsilon ) \frac{C_1}{\mu }  \sinh (\mu t)   \quad & \mbox{\rm for all}&  \quad t \in [a,b),\\
\dsp \dot F(t) & \geq &   \hspace{- 0.2cm} C_1\cosh(\mu t)+ \mu C_0  \sinh (\mu t)  \quad & \mbox{\rm for all} & \quad t \in [a,b),\\
\dsp \ddot F(t)  & \geq & \hspace{-2.4cm} \delta _0   \nu_\phi   (t)^{ -1}   F(t) ^{p}\quad & \mbox{\rm for all} & \quad t \in [a,b),   
 \end{array} \right. 
\] 
with large $a$. The Lemma~\ref{L_ODIEXP} with $A(t)= \cosh(\mu t) $  and $\nu (t)= \cosh(p\mu t)e^{-a_\phi t} t^{-b_\phi } $ shows that if $F(t) \in C^2([0,b)) $, then $b$ must be finite.  
\smallskip

Now consider the case of $\mu =0$.  
Let $ C_0$ and $C_1 $ be defined in (\ref{2.10_Mink}), while the function $\phi _0(x,t) $ satisfies 
\begin{eqnarray*}
\int_{{\mathbb R}^n}  \phi _0(x,t)\,  dx = C_0+C_1 t \, .
\end{eqnarray*}
Then  Corollary~\ref{CLt_tau}   implies 
\begin{eqnarray*}
\int_{{\mathbb R}^n} G[\Gamma (\cdot )|\phi |^{p-1}\phi  ](x,t) \, dx 
& = &
\int_{ 0}^{t}    \Big( \int_{{\mathbb R}^n}  |\phi (z,b )|^{p-1} \phi (z,b ) dz \Big)
(t-b ) \,  db .
\end{eqnarray*} 
Hence,  
\begin{eqnarray*} 
\int_{{\mathbb R}^n}  \phi (x,t)  \,dx 
& =  &
\int_{{\mathbb R}^n}  \phi _0(x,t)  \,dx 
- \int_{ 0}^{t}    \Big( \int_{{\mathbb R}^n}  |\phi (z,b )|^{p-1} \phi (z,b ) dz \Big)
(t-b ) \,  db  \,.  
\end{eqnarray*}
Thus
\begin{eqnarray*}
F(t)
& = &
 C_0+ C_1 t   - \, 
 \int_{ 0}^{t}   \Big( \int_{{\mathbb R}^n}  |\phi (z,b )|^{p-1} \phi (z,b ) dz \Big)
(t-b ) \,  db \,,
\end{eqnarray*}
where $ F(t) $ is defined by (\ref{F}).
It follows $F \in C^2([0,\infty))$. More precisely,
\begin{eqnarray*}
\dot F(t) 
& = &
   C_1    - \, 
 \int_{ 0}^{t}     \Big( \int_{{\mathbb R}^n}  |\phi (z,b )|^{p-1} \phi (z,b ) dz \Big)
 \,  db \,,\\
\ddot F(t) 
& = &
- \, 
 \Big( \int_{{\mathbb R}^n}  |\phi (z,t )|^{p-1} \phi (z,t ) dz \Big)
 \,  db  \,.
\end{eqnarray*}
In particular, for every given positive $\varepsilon <1$ one has
\begin{eqnarray}
\label{48}
F(t)
& \geq &
 C_0+ (1-\varepsilon ) C_1 t \quad \mbox{\rm for all large } \,\, t   \,.
\end{eqnarray}
Indeed, according to (\ref{sign_Mink}) there exists a positive  number   $A_\phi>0 $ such that
\[
\int_{{\mathbb R}^n}  |\phi (x,b) |^{p-1}\phi (x,b)\, dx \leq 0 \qquad 
\mbox{\rm for all} \,\, \, t \geq A_\phi \,.
\]
At the meantime, for every given positive $\varepsilon $ we have 
\begin{eqnarray*}
  \int_{ 0}^{t  }  \left( \int_{{\mathbb R}^n}  |\phi (x,b) |^{p-1}\phi (x,b)\, dx \right) 
 (t-b )  \,  db    
& = &
 \int_{ 0}^{ A_\phi }  dz \int_0^z  \left( \int_{{\mathbb R}^n}  |\phi (x,b) |^{p-1}\phi (x,b)\, dx \right) \,  db \\
&  &
+ \int_{A_\phi }^{t }  dz \int_0^z  \left( \int_{{\mathbb R}^n}  |\phi (x,b) |^{p-1}\phi (x,b)\, dx \right) \,  db \\
& \leq   &
 \int_{ 0}^{ A_\phi }  dz \int_0^z  \left( \int_{{\mathbb R}^n}  |\phi (x,b) |^{p-1}\phi (x,b)\, dx \right) \,  db \\
& \leq   &
\varepsilon t \quad 
\mbox{\rm for large} \,\, \, t \,.
\end{eqnarray*}
The last inequality proves (\ref{48}).
Further, according to  Definition~\ref{D1.2}    
we obtain   
\begin{eqnarray*}
\left| \int_{{\mathbb R}^n} \phi  (x,t) \,dx \right|^{p} 
& \le &
C_\phi  \nu _\phi (t)    \ddot{F}(t)\quad \mbox{\rm for all large} \,\, t  \,,
\end{eqnarray*}
where $  \nu _\phi (t):=   t^{b_\phi }$, $b_\phi \leq 1+p $. Thus
\[
\ddot{F}(t) 
\ge   C_\phi ^{-1}\nu _\phi (t)^{-1} |F(t)|^{p}\quad \mbox{\rm for all large} \,\, t  \,.
\]
By means of the condition \,$C_1   >0$\, we obtain $F(t) >0$ for large $t$ and, 
consequently, 
\[
\ddot{F}(t) 
\ge  \delta _0 \nu _\phi (t)^{-1}   F(t) ^{p}\qquad 
\mbox{\rm for large} \,\,  t \quad \mbox{\rm with} \,\, \delta _0:=  C_\phi ^{-1}>0 \,.
\]
The last inequality together with (\ref{48})
 implies 
 \begin{eqnarray*}
\left\{ \begin{array}{ccccc}
 \dsp F(t)& \geq &   C_0+ (1-\varepsilon )C_1 t \quad &\mbox{\rm for all} &  t \in [a,b),\\
\dsp \ddot F(t) & \geq & \delta _0 \nu _u(t)^{-1}  F(t) ^{p}\quad & \mbox{\rm for all}&  t \in [a,b),    
 \end{array} \right. 
  \end{eqnarray*}
with a sufficiently large number $a$. The  Kato's Lemma~2~\cite{Kato1980} shows that if $F(t) \in C^2([0,b)) $ 
and $\nu _\phi (t)^{-1}   \geq t^{-1-p } $ with $p  >1$,  then $b$ must be finite.
The theorem is proven. \hfill $\square$

\begin{remark}
For the smooth sign preserving function with the support in the ball of radius $r$, 
in the case of $n=3$, one has $\left|\int_{{\mathbb R}^3} \varphi   (x) \, dx  \right|^{3 }\leq   
Cr^6 \left|\int_{{\mathbb R}^3}   \varphi^{3}  (x) \, dx \right| $. 
Hence, if $\phi  $ obeys a finite propagation speed property and its initial values have compact supports, then 
the inequality of  Definition~\ref{D1.2}  is satisfied since
\begin{eqnarray*} 
 \left|\int_{{\mathbb R}^3} \phi  (x,t) \, dx  \right|^{3}
\leq   C_\phi (1+ t)^{6}\left| \int_{{\mathbb R}^3}    \phi^{3} (x,t)  \, dx \right| 
\quad \mbox{for all large } \,\, t   
\,.
\end{eqnarray*} 
For the case of nonlinear wave equation with $\mu = 0$ condition $b_\phi \leq 1+p $ is fulfilled if $5 \leq p$.
In fact, for that  critical case    $p=5$  an existence of the global solution for smooth data with sufficiently small energy is known  
(see, e.g.,   Corollary~6.2~\cite{Shatah}).
\end{remark}

The next conclusion from the theorem stated that the global solution has to change sign.
\begin{corollary}
If the smooth local solution with initial data satisfying (\ref{C0C1_Mink}) with $\sigma =1$ ($\sigma =-1$) 
can be prolonged to the global solution, then the global solution cannot be non-positive (non-negative) for all large
time $t$.
\end{corollary}

\begin{corollary}
\label{C4.4_Min}
Let   $\phi =\phi (x,t) $ be a  continuous global solution of the equation (\ref{1.1_Min_Mink}) with the 
Cauchy data $\phi  (x,0)$, $\phi (x,0) \in C_0^\infty$ satisfying (\ref{C0C1_Mink})  with 
$\sigma = 1$
($\sigma = -1$) and such that   its   self-interaction functional is non-negative 
(non-positive) for all $t$ outside either of the sufficiently  small neighborhood of zero if $\mu >0 $, or inside  of some neighborhood of infinity 
if $\mu =0 $.
Then there exists a sequence $\{ t_k\}_{k=1}^\infty$,    $\lim_{k\to \infty} t_k=\infty$, 
such that the solution has a zero inside of the interior of its support on every  hyperplane $t=t_k$, $k=1,2,\ldots$.
\end{corollary}
 Thus, this global solution is an oscillating in time solution. 
In particular, for the  continuous global solutions  the conditions (\ref{2.10_Mink}),(\ref{sign_Mink}),   and (\ref{C0C1_Mink})
 the sufficient conditions for the creation of the bubbles and their existence in the future. 
 \smallskip

\begin{remark}
The  condition  of the theorem  about the existence of sufficiently small neighborhood of zero means that (\ref{nbdzero})
is fulfilled with some $\varepsilon$, $0 \leq \varepsilon<1$.
\end{remark}

The next theorem generalizes Theorem~\ref{T_Higgs_boson_Mink} by embedding a proper weight 
provided that the time slices of the solution have compact supports. 
The linear Klein-Gordon  equation in the Minkowski spacetime preserves a compactness property of the support 
of the solutions on all time  slices, if it is compact on the initial hyperplane. 
The eigenfunctions of the Laplace operator  give a wide choice for the weight functions. 
In the next theorem they have been used to test global solutions of the equation (\ref{1.1_Min_Mink}).

\begin{theorem}
\label{T_Higgs_boson_Mink_harm}
Let $\phi=  \phi (x,t)\in C ({\mathbb R}^n\times [0,\infty))$ be a weak global solution   of  the
real  field  equation (\ref{1.1_Min_Mink}), which for every given time $t>0$  has a compact support in $x$. Let $\psi =\psi (x)$ be a solution
 of the equation\,
$\Delta \psi = \nu \psi$ in \, $  {\mathbb R}^n$, 
with some number $\nu $ such that $\mu ^2+\nu \geq 0 $. 
Denote the integrals (functionals) of the initial values  
of $\phi $ by
\[
C_{0 \psi } (\phi):= \int_{{\mathbb R}^n} \psi  (x)\phi  (x,0 )   dx,  \quad
C_{1 \psi }(\phi):= \int_{{\mathbb R}^n } \psi  (x)\phi_t (x,0 )   dx \,.
\]
Assume that
the $\psi $-weighted  self-interaction functional 
 $-\lambda\int_{{\mathbb R}^3} \psi  (x) |\phi (x)|^{p-1}\phi (x,t)\, dx$  
satisfies
 \begin{eqnarray}
\label{sign_Mink_harm}
 \int_{{\mathbb R}^n} \psi  (x) |\phi (x,t)|^{p-1}\phi (x,t)\, dx  \leq  0    
\,,
\end{eqnarray} 
for all $t$ either outside of the sufficiently small neighborhood of zero
 if $\mu^2+\nu   >0$, or   inside  of some neighborhood of infinity 
if $\mu^2+\nu  =0$. Assume also  that 
\begin{eqnarray}
\label{C0C1_Mink_harm}
\mu_1 C_{0 \psi }(\phi )+C_{1 \psi }(\phi ) >0\,. 
\end{eqnarray}

Then, the global solution $\phi =\phi (x,t)$ cannot be 
an  asymptotically   time-weighted $-\psi ~L^p$-signed  with the weight  $\nu _\phi = e^{a_{\phi,\psi } t} t^{b_{\phi,\psi } }$, where  if $\mu_1 >0$ then either  $ a_{\phi,\psi } < \mu_1 (p-1)$,  or $ a_{\phi,\psi } = \mu_1 (p-1)$ and 
$b_{\phi,\psi } <-2$, where $\mu_1:=\sqrt{\mu ^2+\nu } $, while $a_{\phi,\psi }=0 $ 
and  $b_{\phi,\psi } \leq 1+p$ if $\mu_1 =0 $. 
\end{theorem}

\smallskip

Thus, according to the theorem, {\sl for the global solutions  for any  eigenfunction  $\psi =\psi (x) $ of the Laplace operator in 
$ {\mathbb R}^n$  the conditions (\ref{sign_Mink_harm}),  (\ref{C0C1_Mink_harm}),  and the inequality of Definition~\ref{D3} cannot hold simultaneously.} 
\smallskip

To check solution for the subject of a zero, one can choose a positive eigenfunction $\psi $ 
of the Laplace operator in ${\mathbb R}^n $ constructed, for example, in 
Lemma~3.1~\cite{Yordanov-Zhang}. Moreover, that eigenfunction $\psi $ has an exponential growth and, consequently, it allows to generalize the last theorem
to the solutions decaying exponentially at infinity.  Then, for an arbitrary positive number $\nu  $ 
existence of the positive   or negative eigenfunction $\psi $  of the Laplace operator can be proved by the scaling arguments.
\smallskip

The proof of  Theorem~\ref{T_Higgs_boson_Mink_harm} is very similar to the one of Theorem~\ref{T_Higgs_boson_Mink} with some modifications based on the following lemma.

\begin{lemma}
\label{L2.6}
  Assume that the smooth function $f=f (x,t)$  for every given time $t>0$ has a compact support. Let $\psi =\psi (x)$ be a solution
 of the equation\, $
\Delta \psi = \nu \psi$\,   in \,${\mathbb R}^n$ 
with some number $\nu \in {\mathbb R}$. Then,  \mbox{\rm (i)} if $\mu ^2+\nu >0 $, then for all $t>0$ we have 
\[
 \int_{{\mathbb R}^n}  \psi (x)G[f](x,t) \, dx
 = 
   \int_{ 0}^{t}  \left( \int_{{\mathbb R}^n}  \psi (x) f (x,b) \, dx \right) 
 \frac{1}{\sqrt{\mu^2+\nu  }} \sinh (\sqrt{\mu^2+\nu  } (t-b ) ) \,  db \,;
\]
(ii) if $\mu ^2+\nu = 0 $, then for all $t>0$ we have
\[ 
 \int_{{\mathbb R}^n}  \psi (x)G[f](x,t) \, dx
 = 
   \int_{ 0}^{t}  \left( \int_{{\mathbb R}^n}  \psi (x) f (x,b) \, dx \right) 
  (t-b )  \,  db\,;
\]
(iii) if $\mu ^2+\nu < 0 $, then for all $t>0$ we have 
\[ 
 \int_{{\mathbb R}^n}  \psi (x)G[f](x,t) \, dx
 = 
   \int_{ 0}^{t}  \left( \int_{{\mathbb R}^n}  \psi (x) f (x,b) \, dx \right) 
 \frac{1}{\sqrt{|\mu^2+\nu|  }} \sin (\sqrt{|\mu^2+\nu|  } (t-b ) ) \,  db\,.
\]
\end{lemma}
\smallskip

\noindent
{\bf Proof of lemma.}
Let $\phi =\phi (x,t)$ be a  function defined as follows:
\[
\phi (x,t):= G[f ](x,t) \,.
\]
For every given time $t>0$ it has a compact support. Let $\psi =\psi (x)$ be a solution of the equation $
\Delta \psi = \nu \psi$\,   in \,${\mathbb R}^n$.  
We integrate the   identity $\psi (x)\phi (x,t) =  \psi (x)G[f](x,t) $ and obtain 
\[
\int_{{\mathbb R}^n} \psi (x) \phi  (x,t)\, dx  = \int_{{\mathbb R}^n}  \psi (x)G[f ](x,t) \, dx  \,.
\]
Thus, for $F_\psi (t):= \int_{{\mathbb R}^n} \psi (x)   \phi  (x,t)    \, dx  $ we have
\[
 F_\psi(t)=  
 \int_{{\mathbb R}^n}  \psi (x)G[f](x,t) \, dx  \,.
\]
The function $\phi =\phi (x,t)$ solves the Cauchy problem 
\[
\phi  _{tt} -  \bigtriangleup \phi    - \mu ^2 \phi  =   f,   \quad \phi    (x,0) =  0 , \quad \partial_{t }\phi    (x,0 ) =0\,.
\]
From the last equation we derive
\[
\frac{d^2}{d t^2} \int_{{\mathbb R}^n} \psi (x) \phi  (x,t)\, dx    -   \int_{{\mathbb R}^n} \psi (x) \Delta  \phi (x,t) \, dx 
-   \mu^2 \int_{{\mathbb R}^n} \psi (x)  \phi  (x,t)  \, dx=     
\int_{{\mathbb R}^n} \psi (x)    f (x,t)  \, dx 
\] 
and, consequently,
\[
\frac{d^2}{d t^2} \int_{{\mathbb R}^n} \psi (x) \phi  (x,t)\, dx     - 
 ( \mu^2 + \nu )\int_{{\mathbb R}^n} \psi (x)   \phi  (x,t)    \, dx   =    \int_{{\mathbb R}^n} \psi (x)  f(x,t)    \, dx \,.
\]
In the case of  $\mu ^2+\nu >0 $  it follows then that the function $F_\psi(t)$ is
\[
F_\psi(t)
 = 
\int_{ 0}^{t}  \left( \int_{{\mathbb R}^n}  \psi (x)   f (x,b)  \, dx \right) 
 \frac{1}{\sqrt{\mu^2+\nu  } } \sinh (\sqrt{\mu^2+\nu  }(t-b ) ) \,  db    \,.
\]
The remaining cases also follow in a similar manner. The lemma is proven. \hfill $\square$
\smallskip

\begin{corollary}  
\label{C2.7}
Assume that the smooth function $f=f (x,t)$  for every given time $t>0$ has a compact support. Let $\psi =\psi (x)$ be a
harmonic function in ${\mathbb R}^n $. Then
\begin{eqnarray*}
 \int_{{\mathbb R}^n}  \psi (x)G[f](x,t) \, dx
&  = &
   \int_{ 0}^{t}  \left( \int_{{\mathbb R}^n}  \psi (x) f (x,b) \, dx \right) 
 \frac{1}{ \mu   }  \sinh ( \mu    (t-b ) ) \,  db , \quad \mu > 0 \,,\\
 \int_{{\mathbb R}^n}  \psi (x)G[f](x,t) \, dx
& = &
  \int_{ 0}^{t}  \left( \int_{{\mathbb R}^n}  \psi (x) f (x,b) \, dx \right) 
  (t-b )   \,  db, \quad \mu =0 ,\\
 \int_{{\mathbb R}^n}  \psi (x)G[f](x,t) \, dx
&  = &
   \int_{ 0}^{t}  \left( \int_{{\mathbb R}^n}  \psi (x) f (x,b) \, dx \right) 
 \frac{1}{ |\mu|   }  \sin  ( |\mu|    (t-b ) ) \,  db , \quad \mu^2 < 0 \,.
\end{eqnarray*}
In particular,
\begin{eqnarray*}
 \int_{{\mathbb R}^n}  G[f](x,t) \, dx
&  = &
   \int_{ 0}^{t}  \left( \int_{{\mathbb R}^n}   f (x,b) \, dx \right) 
 \frac{1}{\mu } \sinh (\mu (t-b ) ) \,  db, \quad \mu > 0 \,,\\
 \int_{{\mathbb R}^n}  G[f](x,t) \, dx
&  = &
  \int_{ 0}^{t}  \left( \int_{{\mathbb R}^n}   f (x,b) \, dx \right) 
  (t-b ) \,  db, \quad \mu =0 \,,\\
 \int_{{\mathbb R}^n}  G[f](x,t) \, dx
&  = &
   \int_{ 0}^{t}  \left( \int_{{\mathbb R}^n}   f (x,b) \, dx \right) 
 \frac{1}{|\mu |} \sin  (|\mu |(t-b ) ) \,  db, \quad \mu^2 < 0  \,.
\end{eqnarray*}
\end{corollary}
\smallskip

\noindent
{\bf Proof of Theorem~\ref{T_Higgs_boson_Mink_harm}.} 
Let $\phi _0=\phi _0(x,t)$ be a solution of the Cauchy problem (\ref{3.5_Mink}).
Then the weak solution $\phi =\phi (x,t)$ of the equation (\ref{1.1_Min_Mink})  that takes initial values $ \phi  (x,0) = \varphi _0(x), \quad \pa_{t }\phi  (x,0) = \varphi _1(x)$,  
solves the integral equation (\ref{2.6new_Mink}).
It follows
\begin{equation}
\label{43}
\psi (x)\phi (x,t) = \psi (x)\phi _0(x,t)- \psi (x)G[\lambda |\phi  |^{p-1}\phi ](x,t)   \,.
\end{equation}
We have
\begin{eqnarray*}
C_{0 \psi }
& := &
\int_{{\mathbb R}^n}  \psi (x)\phi (x,0 ) \,  dx=  \int_{{\mathbb R}^n}  \psi (x)\varphi _0(x) \,  dx,  \\
C_{1 \psi }
& := &
\int_{{\mathbb R}^n }  \psi (x) \partial_t \phi (x,0 ) \,  dx = \int_{{\mathbb R}^n }  \psi (x) \varphi _1(x) \,  dx\,.
\end{eqnarray*}
If $\mu_1 >0$, then it is easily seen that
\[
\int_{{\mathbb R}^n}\psi (x) \phi _0(x,t)  dx = C_{0 \psi } \cosh (\mu _1 t) +C_{1 \psi }\frac{1}{\mu _1 } \sinh(\mu _1  t) 
\qquad \mbox{\rm for all} \quad t \geq 0\, . 
\]
We integrate (\ref{43}) and obtain 
\[
\int_{{\mathbb R}^n} \psi (x) \phi  (x,t)\, dx  = \int_{{\mathbb R}^n} \psi (x) \phi_0  (x,t)\, dx- \int_{{\mathbb R}^n}  \psi (x)G[\lambda |\phi  |^{p-1}\phi ](x,t) \, dx  \,.
\]
Finally, for the function $F_\psi(t):= \int_{{\mathbb R}^n} \psi (x) \phi  (x,t)\, dx$ we obtain 
\[
 F_\psi(t)=  C_{0 \psi } \cosh (\mu _1 t) +C_{1 \psi }\frac{1}{\mu _1 } \sinh(\mu _1  t) 
- \int_{{\mathbb R}^n}  \psi (x)G[\lambda |\phi  |^{p-1}\phi ](x,t) \, dx  \,.
\]
On the other hand, according to Corollary~\ref{C2.7}, we have 
\[
 \int_{{\mathbb R}^n}  \psi (x)G[ |\phi  |^{p-1}\phi ](x,t) \, dx
= 
   \int_{ 0}^{t}  \left( \int_{{\mathbb R}^n}  \psi (x)    |\phi (x,b) |^{p-1}\phi (x,b)  \, dx \right) 
 \frac{1}{\mu _1 } \sinh (\mu _1(t-b ) ) \,  db\,.
\]
 Thus,
 \begin{eqnarray*}
 F_\psi(t)
& = &
 C_{0 \psi } \cosh (\mu_1 t) +C_{1 \psi }\frac{1}{\mu_1 } \sinh(\mu_1  t) \\
&  &
- \int_{ 0}^{t}  \left( \int_{{\mathbb R}^n}  \psi (x)    |\phi (x,b) |^{p-1}\phi (x,b)  \, dx \right) 
 \frac{1}{\mu _1 } \sinh (\mu _1 (t-b ) ) \,  db  \,.
\end{eqnarray*}
The remaining part of the proof is  similar to the proof of Theorem~\ref{T_Higgs_boson_Mink} and we skip it.
\hfill $\square$ 
\smallskip

We can similarly consider the following equation 
\[
\phi _{tt} -  \bigtriangleup \phi   =  \mu ^2 \phi  - \Gamma (t ) 
\left| \int_{{\mb R}^n} |\phi  (y,t )|^{p-1}\phi  (y,t ) dy \right|^\beta |\phi  |^{p-1}\phi\,,  
\]
which contains  the {\it non-local nonlinearity} ({\it non-local self-interaction}).

\section{Equation in the de~Sitter spacetime. Proof of Theorem~\ref{T2.1}}
\label{S_Proofs}

We consider the case of $\sigma =1$ only, since case of  $\sigma =-1$ follows by reflection $\phi \rightarrow  -\phi  $. Let 
$u_0 \in C^1([0,\infty)\times {\mb R}^n )$ be a function with 
\[ 
\int_{{\mathbb R}^n} u_0(x,t)  dx = C_0 \cosh (Mt) +C_1\frac{1}{M} \sinh(M t) 
\qquad \mbox{\rm for all} \quad t \geq 0\, , 
\]
where
\begin{equation}
\label{2.10}
C_0:= \int_{{\mathbb R}^n} u_0(x,0 ) \,  dx,  \quad
C_1:= \int_{{\mathbb R}^n } \partial_t u_0(x,0 ) \,  dx \,.
\end{equation}
Suppose that $u  \in C([0,\infty); L^q({\mb R}^n))$  is a solution to 
(\ref{2.7}) generated by $u_0 $. 
According to the definition of the solution, for every given $T>0$ we have 
\[
G\left[\Gamma (\cdot ) \left| \int_{{\mb R}^n}|u (y,\cdot )|^{p-1}u(y,\cdot ) \, dy \right|^\beta  |u |^{p-1}u\right]  \in C ([0,T]; L^{q }({\mb R}^n))\cap C ^1([0,T]; {\mathcal D}'({\mb R}^n))\,.
\]
Then   $u  \in C([0,\infty); L^1({\mb R}^n))$ and 
we may integrate the equation (\ref{2.7}):
\begin{equation}
\label{1.7a} 
\int_{{\mb R}^n} u(x,t)  \,dx = \int_{{\mb R}^n} u_0(x,t)  \,dx 
- \int_{{\mb R}^n } G\left[\Gamma (\cdot ) \left| \int_{{\mb R}^n}|u (y,\cdot )|^{p-1}u(y,\cdot ) \, dy \right|^\beta  |u |^{p-1}u\right] (x,t)  \,dx   .
\end{equation}
In particular,
\begin{eqnarray*}
& &
 \int_{{\mb R}^n} u (x,0) \,dx= \int_{{\mathbb R}^n} u_0(x,0 ) \,  dx=C_0, 
\quad \int_{{\mb R}^n} u_t (x,0)\,dx = \int_{{\mathbb R}^n } \partial_t u_0(x,0 ) \,  dx=C_1\,.
\end{eqnarray*}
Consider the case of odd $n \geq 3$. The case of even $n$ can be discussed similarly. Then, 
for the smooth function $u=u(x,t) $     we obtain 
 \begin{eqnarray*} 
& &
\int_{{\mathbb R}^n} G\left[\Gamma (\cdot ) \left| \int_{{\mb R}^n}|u (y,\cdot )|^{p-1}u(y,\cdot ) \, dy \right|^\beta  |u |^{p-1}u\right](x,t) \, dx  \\
 & = &
\int_{{\mathbb R}^n} \, dx \,2\int_{ 0}^{t} db
  \int_{ 0}^{ e^{-b}- e^{-t}} dr_1     \Bigg(  \frac{\partial }{\partial r} 
\Big(  \frac{1}{r} \frac{\partial }{\partial r}\Big)^{\frac{n-3}{2} } 
\frac{r^{n-2}}{\omega_{n-1} c_0^{(n)} }  \\
&  &
\times \int_{S^{n-1} } \left[ \Gamma  (b ) \Big| \int_{{\mathbb R}^n}  |u(z,b )|^{p-1} u(z,b )\, dz \Big|^\beta|u (x+ry,b) |^{p-1}u(x+ry,b) \right]\, dS_y  
\Bigg)_{r=r_1}    \\
  &  & 
\times  (4e^{-b-t})^{-M}
\left( (e^{-t}  + e^{-b} )^2 - r_1^2   \right)^{-\frac{1}{2}+M}
F\left(\frac{1}{2}-M,\frac{1}{2}-M;1; 
\frac{ (e^{-b}- e^{-t})^2-r_1^2}
{  (e^{-b}+ e^{-t})^2-r_1^2} \right)\,. 
\end{eqnarray*}
Therefore,
\begin{eqnarray*} 
& &
\int_{{\mathbb R}^n} G\left[\Gamma (\cdot )\left| \int_{{\mb R}^n}|u (y,\cdot )|^{p-1}u(y,\cdot )dy\right|^\beta |u |^{p-1}u \right](x,t) \, dx  \\
& = &
2\int_{ 0}^{t} db
  \int_{ 0}^{ e^{-b}- e^{-t}} dr_1  \Bigg\{  \frac{\partial }{\partial r} 
\Big(  \frac{1}{r} \frac{\partial }{\partial r}\Big)^{\frac{n-3}{2} } 
\frac{r^{n-2}}{\omega_{n-1} c_0^{(n)} }  \\
&  &
\times \int_{S^{n-1} } \Big[ \Gamma  (b ) \Big| \int_{{\mathbb R}^n}  |u(z,b )|^{p-1} u(z,b )\, dz \Big|^\beta
\Big( \int_{{\mathbb R}^n}  |u (x+ry,b) |^{p-1}u(x+ry,b)\, dx \Big) \Big]\, dS_y  
\Bigg\}_{r=r_1}   \nonumber \\
  &  & 
\times  (4e^{-b-t})^{-M}
\left( (e^{-t}  + e^{-b} )^2 - r_1^2   \right)^{-\frac{1}{2}+M}
F\left(\frac{1}{2}-M,\frac{1}{2}-M;1; 
\frac{ (e^{-b}- e^{-t})^2-r_1^2}
{  (e^{-b}+ e^{-t})^2-r_1^2} \right)   
\end{eqnarray*}
implies,
\begin{eqnarray*}  
&  &
\int_{{\mathbb R}^n} G \left[\Gamma (\cdot )\left| \int_{{\mb R}^n}|u (z,\cdot )|^{p-1}u(z,\cdot )\,dz\right|^\beta |u |^{p-1}u \right] (x,t) \, dx   \\
& = &
2\int_{ 0}^{t} db\Big[ \Gamma  (b ) \Big| \int_{{\mathbb R}^n}  |u(z,b )|^{p-1} u(z,b )\, dz 
\Big|^\beta\left( \int_{{\mathbb R}^n}  |u (x,b) |^{p-1}u(x,b)\, dx \right) 
\Big]\\
&  &
\times   \int_{ 0}^{ e^{-b}- e^{-t}} dr_1 \,  \Big(  \frac{\partial }{\partial r} 
\Big(  \frac{1}{r} \frac{\partial }{\partial r}\Big)^{\frac{n-3}{2} } 
\frac{r^{n-2}}{\omega_{n-1} c_0^{(n)} }  
\times \int_{S^{n-1} } \, dS_y  
\Big)_{r=r_1}   \nonumber \\
  &  & 
\times  (4e^{-b-t})^{-M}
\left( (e^{-t}  + e^{-b} )^2 - r_1^2   \right)^{-\frac{1}{2}+M} 
F\left(\frac{1}{2}-M,\frac{1}{2}-M;1; 
\frac{ (e^{-b}- e^{-t})^2-r_1^2}
{  (e^{-b}+ e^{-t})^2-r_1^2} \right)\,.
 \end{eqnarray*}

We discuss the following    two  cases separately:  with  positive curved mass, $M>0$, and vanishing curved mass, $M=0$, respectively. 
In the case of $M>0$ we  apply Proposition~\ref{Lt_tau} to evaluate the last term   and obtain 
 \begin{eqnarray*} 
&  &
\int_{{\mathbb R}^n} G[\Gamma (\cdot )\Big| \int_{{\mathbb R}^n}  |u(z,\cdot  )|^{p-1} u(z,\cdot  )\, dz \Big|^\beta|u |^{p-1}u ](x,t)  dx  \\
  & = &
\int_{ 0}^{t}  \Gamma  (b )  \Big| \int_{{\mathbb R}^n}  |u(z,b )|^{p-1} u(z,b )\, dz \Big|^\beta \left( \int_{{\mathbb R}^n}  |u (z,b) |^{p-1}u(z,b)\, dz \right) 
 \frac{1}{M} \sinh (M(t-b ) ) \,  db .
\end{eqnarray*} 
Hence, (\ref{1.7a}) reads as follows:
\begin{eqnarray*} 
\int_{{\mathbb R}^n} u(x,t)  \,dx 
& = &
\int_{{\mathbb R}^n} u_0(x,t)  \,dx 
-  \int_{ 0}^{t}  \Gamma  (b ) \Big| \int_{{\mathbb R}^n}  |u(z,b )|^{p-1} u(z,b )\, dz \Big|^\beta\\
&  &
\hspace{1cm}\times \left( \int_{{\mathbb R}^n}  |u (z,b) |^{p-1}u(z,b)\, dz \right)
 \frac{1}{M} \sinh (M(t-b ) ) \,  db \, .  
\end{eqnarray*}
Taking into account (\ref{2.9}) and  (\ref{2.10}) we derive
\begin{eqnarray*}
\int_{{\mathbb R}^n} u(x,t)  \,dx 
& = &
\frac{1}{2}\left(C_0+\frac{C_1}{M} \right)e^{M t} +\frac{1}{2}\left(C_0-\frac{C_1}{M} \right) e^{-M t}  \nonumber\\ 
&  &
-  \int_{ 0}^{t}  \Gamma  (b ) \Big| \int_{{\mathbb R}^n}  |u(z,b )|^{p-1} u(z,b )\, dz \Big|^\beta\\
&  &
\hspace{1cm}\times \left( \int_{{\mathbb R}^n}  |u (z,b) |^{p-1}u(z,b)\, dz \right) 
 \frac{1}{M} \sinh (M(t-b ) ) \,  db   \, .
\end{eqnarray*}
Thus,
\begin{eqnarray*}
F(t)
& = &
   C_0 \cosh(Mt)+ \frac{C_1}{M}  \sinh (Mt)   \\
&  &
-  \int_{ 0}^{t}  \Gamma  (b )\Big| \int_{{\mathbb R}^n}  |u(z,b )|^{p-1} u(z,b )\, dz \Big|^\beta \left( \int_{{\mathbb R}^n}  |u (z,b) |^{p-1}u(z,b)\, dz \right) 
 \frac{1}{M} \sinh (M(t-b ) ) \,  db    \,,
\end{eqnarray*}
where
$
F(t)
  :=  
  \int_{{\mathbb R}^n }  u(x,t)   \,dx $. 
It follows $ F \in C^2([0,\infty))$. More precisely, 
\begin{eqnarray}
\label{1.17a}
\dot F(t) 
& = &
   C_1\cosh(Mt)+ MC_0  \sinh (Mt)    \\
&  &
 -  \int_{ 0}^{t}  \Gamma  (b ) \Big| \int_{{\mathbb R}^n}  |u(z,b )|^{p-1} u(z,b )\, dz \Big|^\beta\left( \int_{{\mathbb R}^n}  |u (z,b) |^{p-1}u(z,b)\, dz \right)
\cosh (M(t-b ) ) \,  db   \,, \nonumber \\ 
\ddot F(t) 
& = &
 M^2F(t)   -    \Gamma  (t ) \Big| \int_{{\mathbb R}^n}  |u(z,b )|^{p-1} u(z,b )\, dz \Big|^\beta \int_{{\mathbb R}^n}  |u (z,t) |^{p-1}u(z,t)\, dz    \,.
\end{eqnarray}
In particular, due to  (\ref{17}) and to $ \Gamma  (t )\geq 0$,  there is a positive number $\varepsilon  $,
$\varepsilon <1 $, such that
\begin{eqnarray*}
F(t) 
& \geq  &
  (1-\varepsilon )\left( C_0 \cosh(Mt)+ \frac{C_1}{M}  \sinh (Mt) \right) \quad \mbox{\rm for large }  \quad t\,.
\end{eqnarray*}
Indeed, due to (\ref{17}) there exist positive $\varepsilon <1$ and $\delta_\phi>0 $ such that
\begin{eqnarray*}
 &  &
\left| \int_{ 0}^{\delta_\phi  }  \Gamma  (b )\Big| \int_{{\mathbb R}^n}  |u(z,b )|^{p-1} u(z,b )\, dz \Big|^\beta\left( \int_{{\mathbb R}^n}  |u (z,b) |^{p-1}u(z,b)\, dz \right) 
 \frac{1}{M} \sinh (M(t-b ) ) \,  db  \right| \\
&  &
\leq \varepsilon \cosh(M t) \quad 
\mbox{\rm for large} \,\, \, t \,. 
\end{eqnarray*}
According to the conditions of the theorem, we have 
$M C _0  + C_1   >0$. 
By means of the last inequality  we   conclude that $F(t)\geq 0 $ for large $t$ and, consequently, 
\begin{eqnarray*}
\ddot F(t)  
& \geq &
-  \Gamma  (t ) \Big| \int_{{\mathbb R}^n}  |u(z,b )|^{p-1} u(z,b )\, dz \Big|^\beta  \int_{{\mathbb R}^n}  |u (z,t) |^{p-1}u(z,b)\, dz   \qquad 
\mbox{\rm for large} \,\, \, t  .
\end{eqnarray*}
On the other hand, using the condition of the theorem  we obtain
\begin{eqnarray*}
&  & 
\left| \int_{{\mathbb R}^n} u (x,t) \,dx \right|^{p} \\
& \le &
 C_u \nu_u  (t)
\int_{{\mathbb R}^n}  |u (x,t) |^{p-1}u(x,t)\, dx \\
& \leq  &
C_u \nu_u  (t) \Gamma (t)^{-1/(\beta +1)} 
\left( -\Gamma (t)   \left| \int_{{\mathbb R}^n} |u (x,t)|^{p-1}u (x,t) dx  \right|^{ \beta } \int_{{\mathbb R}^n} |u (x,t)|^{p-1}u (x,t) dx  \right)^{1/(\beta +1)}\\
& = &
C_u  \nu_u  (t) \Gamma (t)^{-1/(\beta +1)}   \Big(\ddot{F}(t)-  M^2F(t) \Big)^{1/(\beta +1)} \\
& \leq  &
C_u \nu_u  (t) \Gamma (t)^{-1/(\beta +1)}   \ddot{F}(t)^{1/(\beta +1)}   \qquad 
\mbox{\rm for large} \,\, \, t  \,.
\end{eqnarray*}
Here we have used the inequality $ \Gamma  (t ) > 0$. Thus, since $\nu_u  (t)>0 $ we obtain 
\[
\ddot{F}(t) 
\geq   C_u^{-(\beta +1)} \nu_u  (t) ^{-(\beta +1)}  \Gamma (t) |F(t)|^{p(\beta +1)}
\quad \mbox{\rm for all large} \quad t \,.
\]
The inequality $F(t) \geq 0 $  allows us to rewrite  this estimate as follows
\[
\ddot{F}(t) 
\ge  \delta _0 \nu_u  (t)^{-\beta -1}  \Gamma (t)  F(t) ^{p(\beta +1)}\qquad 
\mbox{\rm for all large} \,\, \, t \quad \mbox{\rm with} \,\, \,\delta _0 :=C_u^{-(\beta +1)} >0 \,.
\]
Hence,  taking into account   the last inequality  we arrive at 
 the following system of the ordinary differential inequalities 
\[ 
  \left\{ 
\begin{array}{ccccc}
\dsp F(t) & \geq  & \hspace{- 0.2cm}  (1-\varepsilon )C_0 \cosh(Mt)+ (1-\varepsilon )\frac{C_1}{M}  \sinh (Mt)   \quad & \mbox{\rm for all}&  \quad t \in [a,b),\\
\dsp \dot F(t) & \geq &   C_1\cosh(Mt)+ MC_0  \sinh (Mt)  \quad & \mbox{\rm for all} & \quad t \in [a,b),\\
\dsp \ddot F(t)  & \geq & \hspace{-0.7cm} \delta _0   \nu_u  (t)^{-\beta -1}  \Gamma (t)  F(t) ^{p(\beta +1)}\quad & \mbox{\rm for all} & \quad t \in [a,b),   
 \end{array} \right. 
\] 
with large $a$.  Lemma~\ref{L_ODIEXP} shows that if $F(t) \in C^2([0,b)) $, then $b$ must be finite.

Indeed, we apply Lemma~\ref{L_ODIEXP} with $A(t)= e^{Mt}$ and $p$ replaced with  $p(\beta +1)$. More precisely, if we set 
\[
A(t) = e^{Mt},\qquad \gamma (t) =  \nu_u  (t)^{-\beta -1} \Gamma (t)e^{Mp(\beta +1) t},
\]
then the conditions of Lemma~\ref{L_ODIEXP} read as follows: 
\[
p(\beta +1)>1 \quad \mbox{\rm and } \quad  \Gamma_t   (t)   \leq 0 \quad \mbox{\rm for all } \quad t \in [0,\infty).
\]
The last inequality follows from the monotonicity of $\Gamma (t)$. For the increasing function $\Gamma (t)$ in order to apply 
Lemma~\ref{L_ODIEXP} we replace it with the positive constant, which does not affect the above written 
system of the ordinary differential inequalities. 
 By the condition of the theorem, if the global solution $u=u(x,t)$ is 
an  asymptotically   time-weighted $L^p$-non-positive (-non-negative)  with the weight  $\nu_u$  , then there exist $\varepsilon >0$ and $ c>0$ such that 
\begin{eqnarray*}
& &
\Gamma (t) \geq c  \nu_u  (t)^{\beta +1} e^{-M(p(\beta +1)-1 )t} t^{2+\varepsilon}\quad \mbox{\rm for all} \quad t \in [a,\infty), 
\end{eqnarray*} 
that coincides with (\ref{GammaAp}). The case of $M>0$ is proved.
\medskip

Now consider the case of $M=0$.  
Let 
\begin{eqnarray*} 
 \int_{{\mathbb R}^n}  u_0(x,t)  dx = C_0+C_1 t \, .
\end{eqnarray*}
Then  Corollary~\ref{CLt_tau}   allows us to write 
\begin{eqnarray*} 
&  &
\int_{{\mathbb R}^n} G[\Gamma (\cdot )\Big| \int_{{\mathbb R}^n}  |u(z,\cdot  )|^{p-1} u(z,\cdot  ) dz \Big|^{\beta }|u|^{p-1}u ](x,t) \, dx \\
& = &
\int_{ 0}^{t}  \Gamma  (b ) \Big| \int_{{\mathbb R}^n}  |u(z,b )|^{p-1} u(z,b ) dz \Big|^{\beta }  \Big( \int_{{\mathbb R}^n}  |u(z,b )|^{p-1} u(z,b ) dz \Big)
(t-b ) \,  db .
\end{eqnarray*} 
Hence,  (\ref{1.7a}) reads as follows:
\[ 
\int_{{\mathbb R}^n}  u(x,t)  \,dx 
  =   
\int_{{\mathbb R}^n}  u_0(x,t)  \,dx 
- \int_{ 0}^{t}  \Gamma  (b ) \Big| \int_{{\mathbb R}^n}  |u(z,b )|^{p-1} u(z,b ) dz \Big|^{\beta } \Big( \int_{{\mathbb R}^n}  |u(z,b )|^{p-1} u(z,b ) dz \Big)
(t-b ) \,  db  .  
\]
Thus,
\[
F(t)
  =  
 C_0+ C_1 t   - \, 
 \int_{ 0}^{t}  \Gamma  (b ) \Big| \int_{{\mathbb R}^n}  |u(z,b )|^{p-1} u(z,b ) dz \Big|^{\beta }   \Big( \int_{{\mathbb R}^n}  |u(z,b )|^{p-1} u(z,b ) dz \Big)
(t-b ) \,  db  ,
\]
where $ F(t)   :=   \int_{{\mathbb R}^n}  u(x,t)  \,dx $.
It follows $F \in C^2([0,\infty))$. More precisely,
\begin{eqnarray*}
\dot F(t) 
& = &
   C_1    - \, 
 \int_{ 0}^{t}  \Gamma  (b ) \Big| \int_{{\mathbb R}^n}  |u(z,b )|^{p-1} u(z,b ) dz \Big|^{\beta }  \Big( \int_{{\mathbb R}^n}  |u(z,b )|^{p-1} u(z,b ) dz \Big)
 \,  db \\
\ddot F(t) 
& = &
- \, 
\Gamma  (t )  \Big| \int_{{\mathbb R}^n}  |u(z,t )|^{p-1} u(z,t ) dz \Big|^{\beta }  \Big( \int_{{\mathbb R}^n}  |u(z,t )|^{p-1} u(z,t ) dz \Big)
 \,  db  \,.
\end{eqnarray*}
In particular, with some positive  $\varepsilon <1 $ we have
\begin{eqnarray}
\label{3.21}
F(t)
& \geq &
 C_0+ (1-\varepsilon )C_1 t \quad  \mbox{\rm for large} \quad t\,.
\end{eqnarray}
On the other hand, according to the conditions of the theorem, we obtain
\begin{eqnarray*}
\left| \int_{{\mathbb R}^n} u (x,t) \,dx \right|^{p} 
& \le &
 C_u  \nu _u (t)  \Gamma (t)^{-1/(\beta +1)}   \ddot{F}(t)^{1/(\beta +1)}  \,,
\end{eqnarray*}
where $ \nu _u (t)=t^{b_u } $, $b_u \leq 1+p $. Thus
\[
\ddot{F}(t) 
\ge     C_u ^{-(\beta +1)}\nu _u (t)^{-(\beta +1)}\Gamma (t) |F(t)|^{p(\beta +1)}
\]
for all large $t$. By means of the condition \,$C_1   >0$\,
 we   conclude 
\[
\ddot{F}(t) 
\ge  \delta _0 \nu _u(t)^{-(\beta +1)}\Gamma (t)  F(t) ^{p(\beta +1)}\qquad 
\mbox{\rm for large} \,\,  t  
\]
with $\delta _0:=  C_u^{-(\beta +1)}>0$. The last inequality together with (\ref{3.21})
 implies 
 \begin{eqnarray*}
\left\{ \begin{array}{ccccc}
 \dsp F(t)& \geq &  \hspace{-1cm} C_0+ (1-\varepsilon )C_1 t \quad &\mbox{\rm for all} &  t \in [a,b),\\
\dsp \ddot F(t) & \geq & \delta _0 \nu _u(t)^{-(\beta +1)} \Gamma (t)  F(t) ^{p(\beta +1)}\quad & \mbox{\rm for all}&  t \in [a,b),    
 \end{array} \right. 
  \end{eqnarray*}
with some $a$. The  Kato's Lemma~2~\cite{Kato1980} shows that if $F(t) \in C^2([0,b)) $ and $\nu _u(t)^{-(\beta +1)} \Gamma (t) \geq t^{-1-p(\beta +1)} $ with $p(\beta +1) >1$,  then $b$ must be finite.
 The theorem is proven. \hfill $\Box$
\begin{remark}
In fact, we have proved that any solution $u =u(x,t)$ with permanently bounded support blows up if  $MC_0+C_1>0$ and $M\geq 0$.
\end{remark}

\section{The Higgs boson in the de~Sitter spacetime. Proof of Theorem~\ref{T2.1_psi}} 
\label{SHiggs}

To prove Theorem~\ref{T2.1_psi} we have to apply Lemma~\ref{L2.6} with $\mu ^2+\nu  $ replaced with $ M^2+\nu $, and 
follow the outline of the proof of Theorem~\ref{T_Higgs_boson_Mink_harm}. We leave details of the proof to the reader.

For the differential equation (\ref{Higgs_eq}) of the  Higgs boson in the de~Sitter spacetime, Theorem~\ref{T2.1_psi} leads to the following result,
which will be discussed in the remaining part of this paper. 
\begin{theorem}
\label{T_Higgs_boson}
Let $\phi=  \phi (x,t)\in C ({\mathbb R}^3\times [0,\infty))$ be a weak global solution   of  the
real  field  equation (\ref{Higgs_eq}). Let $ \psi =\psi (x)$ be an eigenfunction of the 
Laplace operator in ${\mathbb R}^3 $ corresponding to the eigenvalue $\nu  $.
Denote by  
\[
C_0(\phi,\psi ):= \int_{{\mathbb R}^3} \psi (x) \phi  (x,0 )   dx,  \quad
C_1(\phi, \psi ):= \int_{{\mathbb R}^3 } \psi (x)\phi_t (x,0 )   dx \,,
\] 
the integrals (functionals) of its $\psi $-weighted initial values and assume    that 
\begin{eqnarray}
\label{C0C1} 
\left(\sqrt{9  + 4(\mu ^2+ \nu )}+3 \right) C_0(\phi,\psi  )+2C_1(\phi,\psi  )  >0. 
\end{eqnarray}  
 Assume also that 
the  $\psi $-weighted  self-interaction functional \,\,
 $-\lambda\int_{{\mathbb R}^3}  \psi (x)  \phi ^{3}(x,t)\, dx$ satisfies
 \begin{eqnarray}
\label{sign}
\int_{{\mathbb R}^3}   \psi (x) \phi ^{3}(x,t)\, dx  \leq  0 
\end{eqnarray} 
for all $t$ either  outside of the sufficiently  small neighborhood of zero if $\mu^2+\nu >0 $,
or inside of some neighborhood of infinity if $\mu^2+\nu =0 $. 

Then, the global solution $\phi =\phi (x,t)$ cannot be 
an  asymptotically   time-weighted $-\psi ~L^3$-signed  with the weight  $\nu _\phi = e^{a_{\phi,\psi } t} t^{b_{\phi,\psi } }$, where if $\mu^2+\nu >0 $, then  either  $a_{\phi,\psi } < \sqrt{9+4(\mu^2+\nu )}-3$   
or $ a_{\phi,\psi } = \sqrt{9+4(\mu^2+\nu )}-3$ and 
$b_{\phi,\psi } <-2$, while $a_{\phi,\psi }=0$ and $b_{\phi,\psi } \leq 4 $ 
if $\mu^2+\nu =0 $. 
\end{theorem}
\smallskip

\noindent
{\bf Proof.}  In this case $\Gamma (t)= e^{-3t}$, $M=\sqrt{9/4+\mu ^2} $, and one can apply Theorem~\ref{T2.1}.
\hfill $\square$
\medskip

Note, in (\ref{C0C1}) the constants $C_0(\phi,\psi  ) $ and $C_1(\phi,\psi  ) $ 
can be arbitrarily small. Then, the set of functions $\phi  $ and $\psi  $ with the properties 
(\ref{C0C1}),  (\ref{sign}) and the inequality of Definition~\ref{D3},
is invariant under action of the  multiplicative group of positive numbers, 
and, consequently, is a conic set. In particular, for  every 
given positive number $\varepsilon  $ 
the action $ \psi$ to $\varepsilon \psi $, shows that the constant 
$C_{\phi,\psi } $ can be 
made arbitrarily small. 
Then, the transform $ \psi$ to $-\psi $, changes signs in all  inequalities with the opposite signs. 
\smallskip

\begin{corollary}
For the global solutions  
(\ref{C0C1}),  (\ref{sign}) and the inequality of Definition~\ref{D3}
cannot hold simultaneously. 
\end{corollary}
\smallskip

It must be noted that the range of the number $b_{\phi,\psi }  $ jumps when 
$\nu \rightarrow -\mu^2$. Moreover, in (\ref{sign}) the {\it sufficiently small} neighborhood is changed with {\it some} neighborhood
as $\nu \rightarrow -\mu^2$.
That reveals some kind of resonance phenomena.
\smallskip

According to the theorem there is no global in time non-positive solution 
$\phi =\phi (x,t) $ to the equation (\ref{Higgs_eq}) 
such that $ (\sqrt{9  + 4\mu ^2}+3) C_0(\phi )+2C_1(\phi )  >0 $.
Indeed, in that case, to verify the last statement, we set $\psi (x)\equiv 1 $.
Analogously, there is no global in time non-negative solution $\phi =\phi (x,t) $ 
to the equation (\ref{Higgs_eq}) such that 
$(\sqrt{9  + 4\mu ^2}+3) C_0(\phi )+2C_1(\phi )  <0 $. 
Hence, we have proved the following result.

\begin{corollary}
Let   $\phi =\phi (x,t) $ be a non-trivial local in time solution of the equation 
(\ref{Higgs_eq}) with the Cauchy data $\phi  (x,0)$, $\phi_t (x,0) \in C_0^\infty$ 
satisfying (\ref{C0C1})  with $\psi (x) \equiv1$
($\psi (x) \equiv -1$) and (\ref{sign})
 for all $t$ outside of the sufficiently  small neighborhood of zero.
Then that local solution cannot be prolonged  to the global solution, which is 
non-positive (non-negative) for all large $t$.
\end{corollary}

 Thus,  the continuous global solution  obtained by prolongation of such local solution must change a sign  
and, consequently, it vanishes  at some points. In particular, such radial global solution has zeros and therefore it
gives rise
to  at least one bubble.   
Hence, for the global solutions, Theorem~\ref{T_Higgs_boson} guarantees  the creation of the bubble. 
Moreover, the next
corollary states that the bubbles exist in any neighborhood of infinite time.

\begin{corollary}
\label{C4.4}
Let   $\phi =\phi (x,t) $ be a  continuous global solution of the equation (\ref{Higgs_eq}) with the 
Cauchy data $\phi  (x,0)$, $\phi_t (x,0) \in C_0^\infty$ satisfying (\ref{C0C1})  with 
$\psi (x) \equiv 1$
($\psi (x) \equiv -1$) and such that   its   self-interaction functional is non-negative 
(non-positive) for all $t$ outside of the sufficiently  small neighborhood of zero.
Then there exists a sequence $\{ t_k\}_{k=1}^\infty$,    $\lim_{k\to \infty} t_k=\infty$, 
such that the solution has a zero inside of the interior of its support on every  hyperplane $t=t_k$, $k=1,2,\ldots$.
\end{corollary}
 Thus, the global solution is an oscillating in time solution. 
In particular, for the  continuous global solutions,  the conditions (\ref{sign})  and (\ref{C0C1}) are 
 the sufficient conditions for the creation of the bubbles and their existence in the future. 
 \smallskip

Furthermore, all statements of the above corollaries are also true if the eigenfunction 
$\psi =\psi (x)   $ is non-constant.
 \smallskip

If initial data have compact support, then the support of solution is contained in 
some cylinder $B_R(0)\times [0,\infty) $. Therefore, if  $ \psi (x)\phi (x,t) $ 
does not change sign, then  the inequality of Definition~\ref{D3} is satisfied with 
$a_{\phi,\psi }=
 b_{\phi,\psi }=0$ for  $ \varepsilon \psi (x)\phi (x,t) $, provided that $\varepsilon >0$
is sufficiently small. 
 \smallskip
 
The last case deserves special consideration. Assume that the Cauchy data $\phi  (x,0)$, $\phi_t (x,0) \in C_0^\infty(B_R(0))$.
Then by the finite speed of propagation property for the solution we have \, supp\,$\phi  \subseteq   B_{R+1}(0) \times [0,\infty) $.
Now we choose the function $\psi =\psi (x) $, in particular,  as an eigenfunction of the Laplace operator  in $B_{\wt R}(0)  $, $\wt R \geq R+1 $, with the  Dirichlet 
data 
$\psi (x)|_{|x|=\wt R} =0$. The eigenvalues of such a problem are well-known (see, e.g., \cite{Tikhonov_Sam}):
\[
\nu _{n, k} = - \left( \frac{\rho ^{(n)} _{k}}{{\wt R}}\right)^2 \,\qquad n=0,1,2,\ldots; \quad k=1,2,3,\ldots\,\,,
\]
where $\rho ^{(n)} _{k}$ are the positive zeros of the Bessel function $J_{n+\frac{1}{2}} $, 
that is the positive roots of the equation $J_{n+\frac{1}{2}}(\rho )=0 $. There are $2n+1$ eigenfunctions
belonging to each  eigenvalue $\nu _{n, k}$. 
In fact, $J_{n+\frac{1}{2}} $ can be written via elementary functions (see, e.g., \cite{B-E}).
The corresponding eigenfunctions in the spherical coordinates are
\begin{eqnarray*} 
\psi _{njk}^{1}
& = &
\sqrt{\frac{\pi \wt R}{2\rho ^{(n)} _{k}r}} J_{n+\frac{1}{2}}\left( \rho ^{(n)} _{k}\frac{r}{\wt R} \right) P^{(j)}_{n} (\cos \vartheta )\cos (j\varphi )
\,,
\quad j=0,1,2,\ldots,n;\\
\psi _{njk}^{2}
& = &
\sqrt{\frac{\pi \wt R}{2\rho ^{(n)} _{k}r}}J_{n+\frac{1}{2}}\left( \rho ^{(n)} _{k}\frac{r}{\wt R} \right) P^{(j)}_{n} (\cos \vartheta )\sin (j\varphi )
\,,\quad j=1,2,3,\ldots,n\,\,.
\end{eqnarray*} 
Here $P^{(j)}_{n} (\xi ) $ are the associated Legendre polynomials.   Thus, we arrive at the following theorem.

\begin{theorem}
\label{T_Higgs_boson_Lap}
Let $\phi=  \phi (x,t)\in C ({\mathbb R}^3\times [0,\infty))$, \mbox{\rm supp}\,$\phi  \subseteq   B_{\wt R}(0) \times [0,\infty) $, 
be a weak global solution   of  the
real  field  equation (\ref{Higgs_eq}). Let $ \psi _{njk}^{i} (x)$ be an eigenfunction of the 
Laplace operator with the vanishing  Dirichlet 
data corresponding to the eigenvalue $\nu _{n, k}  $.
Denote by  
\[
C^{\,\,\,i}_{0\,njk} (\phi  ):= \int_{{\mathbb R}^3} \psi _{njk}^{i}  (x) \phi  (x,0 )   dx,  \quad
C^{\,\,\,i}_{1\,njk}(\phi ):= \int_{{\mathbb R}^3 } \psi _{njk}^{i}  (x)\phi_t (x,0 )   dx \,,
\] 
the integrals (functionals) of its $\psi $-weighted initial values and assume  that 
\[ 
\left(\sqrt{9  + 4(\mu ^2+ \nu_{n, k} )}+3 \right) C^{\,\,\,i}_{0\,njk} (\phi,\psi  )+2C^{\,\,\,i}_{1\,njk}(\phi,\psi  )  >0. 
\]  
Assume also that
the  $\psi _{njk}^{i}  $-weighted  self-interaction functional \,\,
 $-\lambda\int_{{\mathbb R}^3}  \psi _{njk}^{i}  (x)  \phi ^{3}(x,t)\, dx$ satisfies
 \begin{equation}
 \label{54}
\int_{{\mathbb R}^3}   \psi _{njk}^{i}  (x) \phi ^{3}(x,t)\, dx  \leq  0 
\end{equation} 
for all $t$ either outside of  the sufficiently  small neighborhood of zero if $\mu^2+\nu_{n, k} >0 $,
or in some neighborhood of infinity if $\mu^2+\nu_{n, k} =0 $. 

Then, the global solution $\phi =\phi (x,t)$ cannot be 
an  asymptotically   time-weighted $-\psi _{njk}^{i} ~L^3$-signed  with the weight  $\nu _\phi = e^{a_{\phi,\psi } t} t^{b_{\phi,\psi } }$, where if\, $\mu^2+\nu_{n, k} >0 $, then  either \, $a_{\phi,\psi } < \sqrt{9+4(\mu^2+\nu_{n, k} )}-3$,  \,
or  \,
$ a_{\phi,\psi } = \sqrt{9+4(\mu^2+\nu_{n, k} )} -3$ \, and  \,
$b_{\phi,\psi } <-2$, while  \,$a_{\phi,\psi }=0$ \, and  \,$b_{\phi,\psi } \leq 4 $ 
if \,$\mu^2+\nu_{n, k} =0 $. 
\end{theorem}

The next corollary describes a resonance case, when $-\mu^2$ coincides with 
some eigenvalue of the Laplace operator with the  Dirichlet 
condition in some ball with a diameter no less than 
the diameter of the spatial trace of the support of the solution.
Although one can always find such eigenvalues, it does not mean that the corresponding conditions of the theorem are satisfied. 

\begin{corollary}
\label{C_Higgs_boson_Lap}
Let $\phi=  \phi (x,t)\in C ({\mathbb R}^3\times [0,\infty))$, \mbox{\rm supp}\,$\phi  \subseteq   B_{\wt R}(0) \times [0,\infty) $, 
be a weak global solution   of  the
real  field  equation (\ref{Higgs_eq}). Let $ \psi _{njk}^{i} (x)$ be an eigenfunction of the 
Laplace operator with the vanishing  Dirichlet 
data corresponding to the eigenvalue $\nu _{n, k}  $. Assume that 
(resonance)
\[
\mu^2=-\nu_{n, k} \,. 
\]
and
\[
3 C_{0\,injk} (\phi  )+C_{1\,injk}(\phi  )  >0. 
\] 
Assume also that
the  $\psi _{njk}^{i}  $-weighted  self-interaction functional \,\,
 $-\lambda\int_{{\mathbb R}^3}  \psi _{njk}^{i}  (x)  \phi ^{3}(x,t)\, dx$ satisfies (\ref{54}) 
for all $t$  in some neighborhood of infinity. 

Then, the global solution $\phi =\phi (x,t)$ cannot be 
an  asymptotically   time-weighted $-\psi _{njk}^{i} ~L^3$-signed  with the weight  $\nu _\phi =  t^{b_{\phi,\psi } }$, where  \,$b_{\phi,\psi } \leq 4 $. 
\end{corollary}

\smallskip

We note here that, for functions with compact support, the H\"older inequality allows us to verify that all 
conditions of   Theorem~\ref{T_Higgs_boson_Lap} and 
Corollary~\ref{C_Higgs_boson_Lap} are fulfilled as long as the function $\psi _{njk}^{i}  (x)\phi  (x,t)  $ preserves its sign.

\section{Appendix} 
\label{Appendix} 
\subsection{Integral representations for the hyperbolic sine function } 
\label{S_sinh}

In \cite[Sec. 2.4]{B-E} one can find one-dimensional integrals involving hypergeometric function. In this section we
present one more example of such an integral as well as examples of multidimensional integrals appearing in the fundamental solutions
for the Klein-Gordon equation in the de~Sitter spacetime. One can find more examples related to the Tricomi and Gellerstedt equations in
\cite{YagTricomi_GE}, \cite{yagdjian_Rend_Trieste}. 
\begin{proposition} \mbox{\rm \cite{Yag_DCDS}}
\label{Lt_tau}
The function $M^{-1}\sinh (M(t-b))$, $M>0$, with $t\ge b \ge 0$, 
can be represented as follows:\\
\mbox{\rm (i)} In the form of a  one-dimensional integral  
\begin{eqnarray*} 
 \frac{1}{M} \sinh (M(t-b ) )   
& = &
  \int_{  - (e^{-b}- e^{-t})}^{e^{-b}- e^{-t}} 
(4e^{-b-t })^{-M} \Big((e^{-t }+e^{-b})^2 - z^2\Big)^{-\frac{1}{2}+M    }   \\
&   &
\times F\Big(\frac{1}{2}-M   ,\frac{1}{2}-M  ;1; 
\frac{ ( e^{-b}-e^{-t })^2 -z^2 }{( e^{-b}+e^{-t })^2 -z^2 } \Big) \,   dz . \nonumber
\end{eqnarray*}
\mbox{\rm (ii)} If $n$ is odd,  $n=2m+1$, $m  \in {\mathbb N}$, then with $c_0^{(n)} =1\cdot 3\cdot \ldots \cdot (n-2 )$, 
\begin{eqnarray*}
& &
 \frac{1}{M} \sinh (M(t-b ) )  \\ 
& = &
2
  \int_{ 0}^{ e^{-b}- e^{-t}} dr_1 \,  \left(  \frac{\partial }{\partial r} 
\Big(  \frac{1}{r} \frac{\partial }{\partial r}\Big)^{\frac{n-3}{2} } 
\frac{r^{n-2}}{\omega_{n-1} c_0^{(n)} }  \!\!\int_{S^{n-1} } \, dS_y  
\right)_{r=r_1}  \!\! (4e^{-b-t})^{-M} \nonumber \\
  &  & 
\times  
\left( (e^{-t}  + e^{-b} )^2 - r_1^2   \right)^{-\frac{1}{2}+M}
F\left(\frac{1}{2}-M,\frac{1}{2}-M;1; 
\frac{ (e^{-b}- e^{-t})^2-r_1^2}
{  (e^{-b}+ e^{-t})^2-r_1^2} \right). 
\end{eqnarray*} 
\mbox{\rm (iii)} If $n$ is even,  $n=2m$, $m  \in {\mathbb N}$, 
then with $c_0^{(n)} =1\cdot 3\cdot \ldots \cdot (n-1 )$, 
\begin{eqnarray*}
& &
 \frac{1}{M} \sinh (M(t-b ) )  \\
& = &
2\int_{ 0}^{ e^{-b}- e^{-t}}  \!\!dr_1  \!\! \left( \frac{\partial }{\partial r} 
\Big( \frac{1}{r} \frac{\partial }{\partial r}\Big)^{\frac{n-2}{2} } 
\frac{2r^{n-1}}{\omega_{n-1} c_0^{(n)} }  \!\!\int_{B_1^{n}(0)} \frac{1 }{\sqrt{1-|y|^2}}  dV_y 
\right)_{r=r_1}  \nonumber \\
  &  & 
\hspace{2.5cm} \times 
 (4e^{-b-t})^{-M} \left( (e^{-t}  + e^{-b} )^2 - r_1^2   \right)^{-\frac{1}{2}+M}\\
  &  & 
\hspace{2.5cm} \times 
F\left(\frac{1}{2}-M,\frac{1}{2}-M;1; 
\frac{ (e^{-b}- e^{-t})^2-r_1^2}
{  (e^{-b}+ e^{-t})^2-r_1^2} \right)\!\! . 
\end{eqnarray*} 
Here 
the constant $\omega_{n-1} $ is the area of the unit sphere $S^{n-1} \subset {\mathbb R}^n$. 
\end{proposition}
\medskip

If we set $b=0$ in the above integrals, then we get integral representations of the function 
$\sinh (Mt)$ depending on the parameter $M>0$.  
By passing to the limit as $M \to 0$ we arrive at the following corollary.
\begin{corollary}  \mbox{\rm \cite{Yag_DCDS}} 
\label{CLt_tau}
The function $t-b$ with $t\ge b \ge 0$, 
can be represented as follows:\\
\mbox{\rm (i)} In  the form of a one-dimensional integral  
\begin{eqnarray*} 
t-b 
& = &
  \int_{  - (e^{-b}- e^{-t})}^{e^{-b}- e^{-t}} 
\Big((e^{-t }+e^{-b})^2 - z^2\Big)^{-\frac{1}{2}}  F\Big(\frac{1}{2} ,\frac{1}{2};1; 
\frac{ ( e^{-b}-e^{-t })^2 -z^2 }{( e^{-b}+e^{-t })^2 -z^2 } \Big) \,   dz .  \nonumber 
\end{eqnarray*}
\mbox{\rm (ii)} If $n$ is odd,  $n=2m+1$, $m  \in {\mathbb N}$, then with $c_0^{(n)} =1\cdot 3\cdot \ldots \cdot (n-2 )$, 
\begin{eqnarray*}
t-b 
& = &
2
  \int_{ 0}^{ e^{-b}- e^{-t}} dr_1 \,  \left(  \frac{\partial }{\partial r} 
\Big(  \frac{1}{r} \frac{\partial }{\partial r}\Big)^{\frac{n-3}{2} } 
\frac{r^{n-2}}{\omega_{n-1} c_0^{(n)} }  \!\!\int_{S^{n-1} } \, dS_y  
\right)_{r=r_1}   \nonumber \\
  &  & 
\times 
\left( (e^{-t}  + e^{-b} )^2 - r_1^2   \right)^{-\frac{1}{2}}
F\left(\frac{1}{2},\frac{1}{2};1; 
\frac{ (e^{-b}- e^{-t})^2-r_1^2}
{  (e^{-b}+ e^{-t})^2-r_1^2} \right). 
\end{eqnarray*} 
\mbox{\rm (iii)} If $n$ is even,  $n=2m$, $m  \in {\mathbb N}$, 
then with $c_0^{(n)} =1\cdot 3\cdot \ldots \cdot (n-1 )$, 
\begin{eqnarray*}
t-b 
& = &
2\int_{ 0}^{ e^{-b}- e^{-t}} dr_1 \,  \left( \frac{\partial }{\partial r} 
\Big( \frac{1}{r} \frac{\partial }{\partial r}\Big)^{\frac{n-2}{2} } 
\frac{2r^{n-1}}{\omega_{n-1} c_0^{(n)} }  \!\!\int_{B_1^{n}(0)} \frac{1 }{\sqrt{1-|y|^2}} \, dV_y 
\right)_{r=r_1}     \\
  &  & 
\times  
\left( (e^{-t}  + e^{-b} )^2 - r_1^2   \right)^{-\frac{1}{2}}
F\left(\frac{1}{2},\frac{1}{2};1; 
\frac{ (e^{-b}- e^{-t})^2-r_1^2}
{  (e^{-b}+ e^{-t})^2-r_1^2} \right)\!\! .  
\end{eqnarray*} 
\end{corollary}
\medskip

\subsection{Second order differential inequalities} 
\label{S_diffineq}

The second order differential inequalities with  power decreasing kernel play a key role in   proving
 blow-up of the solutions to  semilinear equations. Kato's lemma \cite{Kato1980} 
allows us to derive from the inequality 
\[
\ddot w \geq b t^{-1-p} w^p,\qquad p>1, \,\, b>0, \quad t \, \,\,\mbox{\rm large},
\]
the boundedness  of the life-span
of a solution with the property $\dot w  \geq a >0$.   For the equation in the de~Sitter spacetime the kernel 
$ e^{-Mt}$ of the
corresponding ordinary differential inequality decreases exponentially: 
\[
\ddot w \geq b e^{-Mt} w^p,\qquad p>1, \,\, b>0, \,M>0,\quad t \,\,\, \mbox{\rm large}.
\]
There is  a non-trivial global  solution to the last differential inequality. Hence, to  generalize Kato's lemma  
we need  proper supplementary conditions on  the 
involved functions.

\begin{lemma}  \mbox{\rm \cite{Yag_DCDS}} 
\label{L_ODIEXP}
Suppose $F(t) \in C^2([a,b))$, and 
\begin{equation}
\label{2.25c} 
F(t) \ge c_0 A(t) , \quad  \dot F(t) \geq 0 , \quad 
\ddot F(t) \ge \gamma  (t)A(t)^{- p } F(t)^p \quad \mbox{\rm for all} \,\, t \in [a,b) , 
\end{equation}
where  $A, \gamma  \in C^1([a,\infty))$ are non-negative  functions and $p>1$, $c_0>0$. Assume that  
\[
\lim_{t \to \infty} A(t) = \infty\,,
\]
  and that 
\[ 
\frac{d}{dt} \left( \gamma  (t)A(t)^{- p }  \right) \leq 0\quad \mbox{\rm for all} \quad t \in [a,b)\,.
\]
If there exist $\varepsilon >0$ and $ c>0$ such that 
\begin{eqnarray}
\label{GammaAp}
& &
 \gamma  (t) \geq c A(t) (\ln A(t))^{2+\varepsilon}\quad \mbox{\rm for all} \quad t \in [a,b), 
\end{eqnarray}
then $b$ must be finite.
\end{lemma} 
\smallskip

We note here that the equation
\beqst 
\ddot F(t) =   e^{-dt} F(t)^p\,,\quad d>0,   
\eeqst
has a global solution $F(t)=  c_F e^{\frac{d}{p-1}t}$, where $c_F= \left( {d}/(p-1) \right)^{2/(p-1)}$, while the    
corresponding $A(t)=c_A e^{at}$, $a>0$,
and $\gamma (t)= c_\gamma  e^{(pa-d) t}$.  The condition (\ref{GammaAp}) implies   $a> d/(p-1)$. 
On the other hand, the first inequality of (\ref{2.25c}) holds only if  $a\leq  d/(p-1)$. 
\bigskip

\noindent
{\bf Acknowledgments}

\noindent
This paper was completed during my stay at the Technical University   Bergakademie  Freiberg.  I am grateful to Michael Reissig for the invitation to 
Freiberg and for the warm hospitality. 
I express my gratitude   to the Deutsche Forschungsgemeinschaft for   the financial support under grant GZ: RE 961/16-1 
 AOBJ: 577386.

\end{document}